\documentclass[onecolumn,draftcls,12pt]{IEEEtran}

\usepackage[cmex10]{amsmath}
\usepackage{amsfonts,mathtools}
\usepackage{algorithmic}
\usepackage[]{algorithm}
\usepackage{array}
\usepackage{url}
\usepackage{bm}
\usepackage{graphicx}
\usepackage{mathrsfs}
\usepackage{amsthm} 
\usepackage{amsmath,bm}
\usepackage{subfigure}
\usepackage{cite}
\usepackage{color}
\usepackage{float}
\usepackage[font={small}]{caption}
\usepackage{epstopdf}

\usepackage{hyperref}
\usepackage{cleveref}

\makeatletter
\newlength \figwidth
\if@twocolumn
\else
\fi
\makeatother

\theoremstyle{plain}

\newtheorem{lem}{Lemma}
\newtheorem{prop}{Proposition}

\theoremstyle{definition}

\theoremstyle{remark}
\newtheorem{rem}{\bf Remark}

\providecommand{\abs}[1]{\lvert#1\rvert}												
\renewcommand{\b}[1]{\ensuremath{\mathbf{#1}}}		 							
\newcommand{\bs}[1]{\ensuremath{\boldsymbol{#1}}}		 						
\newcommand{\Ex}[1]{\ensuremath{\mathbb{E}[#1]}} 		  
\newcommand{\norm}[1]{\ensuremath{\left\|#1\right\|}}						

\providecommand{\vect}[1]{\text{vec}\left(#1\right)}
\providecommand{\tr}[1]{\text{tr}\left(#1\right)}

\def \Rn {{\mathbb{R}}}
\def \x {{\b{x}}}
\def \y {{\b{y}}}
\def \xc {{\check{\b{x}}}}
\def \X {{\b{X}}}

\def \Xt {{\tilde{\b{X}}}}
\def \Xh {{\hat{\b{X}}}}
\def \L {{\b{L}}}

\def \I {{\b{I}}}
\def \B {{\b{B}}}
\def \Bt {{\b{B}_t^\epsilon}}

\def \Ba {{\b{B}^{\text{a}}}}
\def \C {{\b{C}}}

\def \J {{\b{J}}}
\def \K {{\b{K}}}
\def \Dl {{\bs{\Delta}}}
\def\D {{\b{D}}}
\def \Xc {{\check{\b{X}}}}
\def \A {{\b{A}}}
\def \Y {{\b{Y}}}

\def \Q {{\b{Q}}}
\def \O {{\mathcal{O}}}
\def \E {{\mathcal{E}}}
\def \G {{\mathcal{G}}}
\def \cjt {{\mathcal{C}^j_t}}
\def \bG {{\boldsymbol{\Gamma}}}

\def \N {{\mathcal{N}}}
\def \Cs {{\mathcal{C}}}

\begin{document}

\title{Stochastic Multidimensional Scaling}
\author{Ketan Rajawat,~\IEEEmembership{Member,~IEEE}, and Sandeep Kumar,~\IEEEmembership{Student Member,~IEEE}
\thanks{ Ketan Rajawat and Sandeep Kumar are with the Department of Electrical Engineering, Indian Institute of Technology Kanpur, Kanpur, UP 208016, India, email: \texttt\{ketan, sandkr\}@iitk.ac.in}.}
\maketitle

\vspace{-0.8cm}



\begin{abstract}
Multidimensional scaling (MDS) is a popular dimensionality reduction techniques that has been widely used for network visualization and cooperative localization. However, the traditional stress minimization formulation of MDS necessitates the use of batch optimization algorithms that are not scalable to large-sized problems. This paper considers an alternative stochastic stress minimization framework that is amenable to incremental and distributed solutions. A novel linear-complexity stochastic optimization algorithm is proposed that is provably convergent and simple to implement. The applicability of the proposed algorithm to localization and visualization tasks is also expounded. Extensive tests on synthetic and real datasets demonstrate the efficacy of the proposed algorithm. 
\end{abstract}


\section{Introduction}\label{intro}
Multidimensional scaling addresses the problem of embedding relational data onto a low-dimensional subspace. Originally proposed in the context 
of psychometrics and marketing \cite{mds-book}, MDS and its variants have since found 
applications in social networks 
\cite{maaten2008visualizing,platzer2013visualization,Tang:2015,vanderMaaten2012,agrafiotis2003stochastic}, genomics \cite{tzeng2008multidimensional}, computational chemistry \cite{choi2011browsing}, machine learning\cite{beatty1997dimensionality}, and 
wireless networks \cite{hero-dwmds}.  As an exploratory technique, MDS is often 
used as a first step towards uncovering the structure inherent to 
high-dimensional data. In the context of machine learning and data 
mining, the pairwise dissimilarities are calculated using high- or 
infinite-dimensional nodal attributes, and MDS yields a 
distance-preserving, low-dimensional embedding. Of particular importance 
are the embeddings obtained in two or three dimensional euclidean 
spaces, that serve as perceptual maps for visualizing relationships 
between objects. In the context of social networks, such representations 
reveal interconnections between people and communities, and are often 
more insightful than simpler metrics such as centrality and density. 
Different from the classical MDS framework that utilizes principal 
component analysis, modern MDS formulations are based on the 
minimization of a non-convex stress function \cite{mds-book}. Since the 
stress function is a weighted sum of squared fitting errors, it allows 
for the possibility of missing and noisy dissimilarities. Consequently 
variants of the stress minimization problem have been developed for 
robust MDS \cite{forero2012sparsity}, visualization of time-varying data \cite{xu2013regularized}, and 
cooperative localization of static \cite{kumar2016cooperative,simonetto2014distributed,hero-dwmds,7518679} and mobile 
networks \cite{kumar2016cooperative}. Popular algorithms for solving the stress minimization 
problem include `scaling by majorizing a complicated function' (SMACOF) 
\cite{mds-book}, semidefinite programming \cite{biswas2006semidefinite}, alternating directions 
method of multipliers \cite{7518679,simonetto2014distributed}, and distributed SMACOF \cite{hero-dwmds}.

The attractiveness of the MDS framework has however started to diminish with the advent of the data deluge. Specifically, when embedding $N$ 
objects, the per-iteration complexity and memory requirements of the 
aforementioned algorithms increase at least as $\mathcal{O}(N^2)$, 
making them impractical for large-scale problems. To this end, 
approximate versions of SMACOF have been proposed for large-scale 
visualization applications \cite{baea2012visualization,ingram2009glimmer}. 
Nevertheless, most approximate MDS algorithms are still too complex for 
large-scale data, and cannot be generalized to other applications such 
as cooperative localization of large networks. 

Visualization or 
localization of time-varying data is even more challenging since the 
iterative majorization algorithm must converge at every time 
instant\cite{hero-dwmds,kumar2016cooperative, xu2013regularized}. In 
mobile sensor networks, carrying out a large number of iterations at 
each time instant incurs a tremendous communication overhead, and is 
generally impractical. For instance, the distributed weighted MDS 
approach \cite{hero-dwmds} still requires at least $N$ operations per 
iteration per time instant, which is prohibitive for large networks. For 
large-scale applications, where localization or visualization is 
constrained by the per-iteration complexity and memory requirements, it 
is instead desirable to have an online algorithm. Towards this end, the 
goal is to obtain an adaptive algorithm that processes dissimilarity 
measurements in a sequential or online manner. For instance, an adaptive 
algorithm can allow visualization of large networks by reading and 
processing the pairwise dissimilarities in small batches. Similarly, the 
communication cost required for large-scale network localization can be 
reduced by processing only a few range measurements at a time.

This paper considers the stress minimization problem in a stochastic 
setting, where the dissimilarity measurements and the weights are 
modeled as random time-varying quantities with unknown distributions. 
The first contribution of this paper is a novel stochastic SMACOF 
algorithm that processes the dissimilarities in an online fashion, and 
is therefore applicable to both static and time-varying scenarios (Sec. 
\ref{online-smacof}). The proposed algorithm is not only scalable, but 
is also amenable to a distributed and asynchronous implementation in ad 
hoc networks (Sec. \ref{imp}). As our second contribution, it is shown 
that the trajectory of the stochastic SMACOF algorithm remains close to 
that of an averaged algorithm, which itself converges to a stationary 
point of the stochastic stress minimization problem (Sec. \ref{conv}). 
The analysis borrows tools from spectral graph theory, stochastic 
approximation, and convergence analysis of the SMACOF algorithm. 
Finally, as the third contribution, the performance of the proposed 
algorithm is tested extensively on various synthetic and real-world data 
sets (Sec. \ref{sim}). The numerical tests confirm the applicability of 
the stochastic SMACOF algorithm to a variety of scenarios.

The notation used in this paper is as follows. Bold upper (lower) case 
letters denote matrices (vectors). The ($m,n$)-th entry of a matrix  
$\A$ is denoted by $[\A]_{mn}$. $\I_N$ is the $N \times N$ identity 
matrix, $\mathbf{0}$ denotes the all-zero matrix or vector, and 
$\mathbf{1}$ denotes the all-one matrix or vector, depending on the 
context. For a vector $\x$, $\norm{\x}$ denotes its $\ell_2$ norm. For a 
matrix $\A$, $\norm{\A}$ denotes its Frobenious norm, $\norm{\A}_2$ 
denotes the $\ell_2$ norm, $\text{tr}(\A)$ denotes its trace, and 
$\text{det}(\A)$ denotes its determinant.

\section{Background and Problem Statement}\label{background}

\subsection{Classical MDS and SMACOF}
The classical MDS framework seeks $P$-dimensional embedding vectors $\{\x_n\}_{n=1}^N$, given the pairwise distances or dissimilarities $\{\delta_{mn}\}_{(m,n)\in\E}$, where $\E \subseteq \{(m,n) \mid 1\leq m < n \leq N\}$, between $N$ different nodes or objects, denoted by the set $\N:=\{1,\ldots,N\}$. The embedding vectors, collected into the rows of $\X \in \Rn^{N \times P}$, are estimated by solving the following non-convex optimization problem \cite{mds-book}
\begin{align}\label{staticmds}
\hat{\X} = \arg \min_{\X} \sum_{1\leq m < n \leq N} w_{mn}\left(\delta_{mn}-\norm{\x_m-\x_n}_2\right)^2
\end{align}
where $w_{mn}$ is the weight associated with the measurement $\delta_{mn}$, and is set to zero for all $(m,n) \notin \E$. The non-zero weights can be chosen in a number of ways, depending on the application, and are often simply set to one. The objective function in \eqref{staticmds} is referred to as the stress function, and is henceforth denoted by $\sigma(\X)$. It can be seen that the optimum $\hat{\X}$ obtained in \eqref{staticmds} is not unique, and exhibits translational, rotational, and reflectional ambiguity. 

The stress-minimization problem in \eqref{staticmds} is non-convex, and can be solved up to a local optimum using the well known SMACOF algorithm. Expanding the stress function, we obtain 
\begin{align}
\sigma(\X) &= \sum_{m < n} w_{mn} \left(\delta_{mn}^2 + \norm{\x_m-\x_n}^2 - 2\delta_{mn}\norm{\x_m-\x_n}\right) \\
&=  \sum_{m<n} w_{mn} \delta_{mn}^2 + \text{tr}(\X^T\L\X) - 2\text{tr}(\X^T\B(\X)\X) \label{lossfun}
\end{align}
where, 
\begin{align}
[\L]_{mn} &= \begin{cases} -w_{mn} &  m\neq n \\
\sum_{k=1}^m w_{mk} & m=n \end{cases}\label{lx}\\
[\B(\X)]_{mn} &=  \begin{cases} -\frac{w_{mn}\delta_{mn}}{\norm{\x_m-\x_n}} &  m\neq n, \x_m \neq \x_n \\
0 &  m\neq n, \x_m=\x_n 	\\
-\sum_{k=1}^m [\B(\X)]_{mk} & m=n \end{cases} \label{bx}
\end{align}
The SMACOF algorithm works by iteratively majorizing the last term in \eqref{lossfun} with a linear function and subsequently minimizing the majorized stress function with respect to $\X$. Starting with an initial $\Xh^{(0)}$, the SMACOF update at the $k$-th iteration entails carrying out the following update:
\begin{align}
\Xh^{(k+1)} &= \arg \min_{\X} \text{tr}(\X^T\L\X) - 2\text{tr}(\X^T\B(\Xh^{(k)})\Xh^{(k)}) \label{smacofopt} \\
&=\L^{\dagger}\B(\Xh^{(k)})\Xh^{(k)} \label{smacof}
\end{align}
where \eqref{smacof} follows since $\B(\X)\X$ lies in the range space of $\L$. Observe that since $\L$ is rank-deficient, the solution to \eqref{smacofopt} is not unique. However, when the weights $\{w_{mn}\}$ specify a fully connected graph $\G := (\{1,\ldots,N\}, \E)$, both $\L$ and $\B(\X)$ have rank $N-1$, with the null space of $\L$ being $\mathbf{1}$. Therefore, any solution to \eqref{smacofopt} is of the form $\L^{\dagger}\B(\Xh^{(k)})\Xh^{(k)} + \mathbf{1}c$ for $c \in \Rn$. Further, if the initial $\X^{(0)}$ is chosen such that it is centered at the origin, i.e., $\b{1}^T\X^{(0)} = \b{0}$, the updates in \eqref{smacof} ensure that $\b{1}^T\X^{(k)} = \b{0}$ for all $k \geq 1$. 

\subsection{Stochastic MDS}
This paper considers the MDS problem in a stochastic setting, where the weights, and dissimilarities or distance measurements are random variables with unknown distributions. Specifically, given $\{\delta_{mn}(t)\}$ and $\{w_{mn}(t)\}$, the stochastic stress minimization problem is formulated as
\begin{align}\label{ss}
\min_{\X} \bar{\sigma}(\X) := \sum_{m < n} \Ex{w_{mn}(t)(\delta_{mn}(t)-\norm{\x_m-\x_n})^2}.
\end{align} 
In the absence of the distribution information, the expression for $\bar{\sigma}(X)$ cannot be evaluated in closed-form, and the SMACOF algorithm cannot be applied. Instead, \eqref{ss} must be solved using a stochastic optimization algorithm. Of particular interest are the so-called \emph{online} algorithms that can process the observations $\{\delta_{mn}(t)\},\{w_{mn}(t)\}$ in an incremental manner. Within this context, efficient implementations of the stochastic (sub-)gradient descent (SGD) method have been used to solve very large-scale problems \cite{Bottou}. The SGD updates utilize the subgradient of the instantaneous objective function, and for the present case, take the form:
\begin{align}\label{sgd}
\check{\X}_{t+1} = \check{\X}_t + \mu \left(\B_t(\check{\X}_t) \check{\X}_t - \L_t\check{\X}_t\right)
\end{align}
where $\mu \in (0,1)$ is the learning rate or step size parameter. While the performance of the SGD has been well-studied for convex problems, the same is not true for non-convex problems, such as the one in \eqref{ss}. Indeed, the standard SGD algorithm does not necessarily converge for many non-convex problems \cite{sa2015global}. In the present case also, the SGD method exhibits divergent behavior; see Sec. \ref{sim}. The general-purpose stochastic majorization-minimization method \cite{mairal2013stochastic} is also not applicable in the present case since it requires a strongly convex surrogate function. On the other hand, problem-specific stochastic algorithms have been developed and applied with great success. Examples include the online expectation-maximization and the online matrix factorization approaches \cite{cappe2009line,Bottou}. Along similar lines, the next section details the stochastic version of the SMACOF algorithm, and studies its asymptotic properties.

\section{Online Embedding via Stochastic SMACOF}\label{online-smacof}
\subsection{Algorithm outline}
Given $\{\delta_{mn}(t)\}$ and $\{w_{mn}(t)\}$, and starting with an arbitrary origin-centered $\Xh_0$, the updates for the proposed stochastic SMACOF algorithm take the form,
\begin{align}\label{smdseq}
\Xh_{t+1} &= (1-\mu)\Xh_t + \mu \L_t^{\dagger}\B^{\epsilon}_t(\Xh_t)\Xh_t & t\geq 0
\end{align}

\begin{align}\label{bepsilon}
\text{where,}\quad[\Bt(\X)]_{mn} &=  \begin{cases} -\frac{w_{mn}(t)\delta_{mn}(t)}{\sqrt{\norm{\x_m-\x_n}^2+\epsilon_x}} &  m\neq n \\
	-\sum_{k=1}^N [\Bt(\X)]_{mk} & m=n \end{cases} 
\end{align}
with $\epsilon_x$ being a small positive constant that ensures that the entries of $\Bt(\X)$ stay bounded for all $\X$. The update rule can be viewed a stochastic version of the SMACOF algorithm with the following modifications (a) at each time instant, only one iteration of SMACOF is executed using the modified definition of $\B^{\epsilon}(\X)$ in \eqref{bepsilon}; (b) the estimated coordinates $\Xh_t$ at time $t$ are used for initialization at $t+1$; and (c) the estimated coordinates $\Xh_{t+1}$ are constructed by taking a convex combination of $\hat{\X}_t$ and the SMACOF output. The last modification endows the algorithm with tracking capabilities since the parameter $\mu$ may be interpreted as the forgetting factor, and can be tuned in accordance with the rate of change of $\{\delta_{mn}(t)\}$ and $\{w_{mn}(t)\}$. For example, the embedding at time $t+1$ can be forced to be close to those at time $t$ by setting $\mu \ll 1$. Finally, the proposed update rule subsumes the SMACOF algorithm for static scenarios, where we set $\delta_{mn}(t) = \delta_{mn}$ and $w_{mn}(t) = w_{mn}$ for all $t$, and $\mu = 1$. 

The update rule in \eqref{smdseq} is valid only if the graph $\G_t$ defined by $\{w_{mn}(t)\}$ is connected for all $t \geq 1$. In the case when $\G_t$ has more than one connected component, the coordinates within each component must be updated separately. Let $\cjt$ be the set of nodes belonging to the $j$-th component and $\I^j_t$ be the $\abs{\cjt} \times N$ selection matrix containing the rows of $\I_N$ corresponding to the elements in $\cjt$. Defining $\L_t^{(j)} := \I^j_t\L_t{\I^j_t}^T$ and $\B_t^{\epsilon}(\X_t^{(j)}):= \I^j_t\B_t^{\epsilon}(\X_t){\I^j_t}^T$, the update rule for the nodes in $\cjt$ is given by 
\begin{align}\label{smdseqc}
\Xh^{(j)}_{t+1} = (\I-\mu \J_t)\Xh^{(j)}_t + \mu (\L_t^{(j)})^{\dagger}\B^{\epsilon}_t(\Xh^{(j)}_t)\Xh^{(j)}_t 
\end{align}
where $\J_t := \I-\bm{11}^T/\abs{\mathcal{C}_j(t)}$, is the $\abs{\mathcal{C}_j(t)} \times \abs{\mathcal{C}_j(t)}$ centering matrix which ensures that the coordinate center of each component does not change after the update, i.e., $\bm{1}^T\Xh^{(j)}_{t+1} = \bm{1}^T\Xh^{(j)}_t$. The general update rule 
\begin{align}\label{smdsg}
\Xh_{t+1} = (\I-\mu \L_t^{\dagger}\L_t)\Xh_t + \mu\L_t^\dagger\B_t(\Xh_t)\Xh_t
\end{align}
subsumes the forms specified in \eqref{smdseq} and \eqref{smdseqc}, irrespective of the number of connected components in $\G_t$, since it holds that
\begin{align}
[\L_t\L_t^{\dagger}]_{mn} = \begin{cases} 1-1/\abs{\mathcal{C}_t^j} &m=n \in \mathcal{C}_t^j \\
-1/\abs{\mathcal{C}_t^j} &m\neq n, m,n \in \mathcal{C}_t^j \\
0 & \text{otherwise.} 
\end{cases}
\end{align}

In contrast to the classical SMACOF algorithm, the proposed algorithm is flexible enough to be used in a number of different scenarios. As already discussed, a specific choice of parameters allows us to interpret the SMACOF algorithm as a special case of the proposed algorithm. On the other hand, the stochastic SMACOF can also be used to solve very large-scale MDS problems, where the full set of measurements $\{\delta_{mn}\}$ cannot be processed simultaneously. Instead, it is possible to apply \eqref{smdsg} on a small subset of observations, corresponding to a subgraph $\G_t$. A special case occurs when exactly one edge is chosen per time instant and per cluster, i.e., $\abs{\mathcal{C}_t^j} = 2$, and the updates in \eqref{smdsg} reduce to those in encountered in the stochastic proximity embedding (SPE) algorithm \cite{agrafiotis2003stochastic}, 
\begin{align}
\x_i(t+1) &= (1-\mu)\x_i(t) + \mu\frac{\delta_{ij}(t)}{\norm{\x_i(t)-\x_j(t)}}\x_i(t) \nonumber\\
&+ \mu\left(1-\frac{\delta_{ij}(t)}{\norm{\x_i(t)-\x_j(t)}}\right)\x_j(t) 
\end{align}
and likewise for node $j$. The proposed stochastic SMACOF is therefore a generalization of SPE, applied to components of arbitrary sizes. Since the updates in \eqref{smdsg} for any two clusters $\mathcal{C}_t^j$ and $\mathcal{C}_t^k$ do not depend on each other, the proposed algorithm can also be implemented in a distributed and asynchronous manner. Such an implementation is particularly suited to the range-based localization problems that arise in wireless networks. 

Finally, akin to the classical adaptive filtering algorithms such as LMS, the proposed algorithm can also be applied to time-varying scenarios, i.e., when $\delta_{mn}(t)$ is non-stationary. The applications of interest include localization of time-varying networks, and visualization of time-varying data. In both cases, the first term $(\I-\mu \J_t)\Xh^{(j)}_t$ in the update \eqref{smdsg} serves as a momentum term. That is, a small $\mu$ encourages $\Xh_{t+1}$ to stay close to $\Xh_t$, resulting in a smooth trajectory of $\{\Xh_t\}$. On the other hand, a large value of $\mu$ enables tracking in highly time-varying scenarios, while making the updates sensitive to noise \cite[Ch-21]{sayed2011adaptive}\cite[Ch-9]{solo1994adaptive}.  Further implementation details pertaining to the localization and visualization problems are discussed in Sec. \ref{imp}. Before proceeding with the asymptotic analysis, the following remark is due. 

\begin{rem}
Building further on the link with adaptive algorithms, $\mu$ may be interpreted as a forgetting factor that downweights the past information. When $\mu$ is a constant that is strictly greater than zero, the algorithm forgets the old data exponentially quickly, thus offering superior tracking capability. In contrast, it is possible to have a long-memory version of the algorithm with a time-varying $\mu_t \rightarrow 0$. As $t \rightarrow \infty$, such an algorithm would no longer track the changes in $\delta_{mn}(t)$, and can be applied to a static scenarios where the algorithm can stop once the embeddings converge. While the bounds developed here apply only to the case of constant $\mu >0$, diminishing step size is in fact utilized in Sec. \ref{sim}.
\end{rem}

\subsection{Asymptotic Performance}\label{conv}
In general, establishing convergence of stochastic algorithms for non-convex problems is quite challenging \cite{sa2015global}. Here, the asymptotic performance of the proposed algorithm is established in two steps. First, it is shown that the trajectory of the stochastic SMACOF algorithm stays close to that of an averaged algorithm, in an almost sure sense. This part involves establishing a hovering theorem, and utilizes techniques from stochastic approximation \cite{solo1994adaptive,kushner2003stochastic,borkarstochastic}. Next, it is shown that the averaged algorithm converges to a stationary point of \eqref{ss}. 

\subsubsection{Assumptions}
For the purposes of establishing convergence, a simplified setting is considered, wherein the graph $\G_t$ at each $t$ consists of $N/p \geq 1$ components of size $p$ each. Let $j_m(t):=\{j \mid m \in \mathcal{C}_t^j\}$ be the index of the component to which node $m$ belongs at time $t$, and define $\boldsymbol{\Theta}_t \in \Rn^{N \times N}$ such that
\begin{align}
[\boldsymbol{\Theta}_t]_{mn} := \begin{cases} -1/N &j_m(t) \neq j_n(t) \\
-1/N+\mu/p & \hspace{-1cm}j_m(t) = j_n(t), m\neq n \\
(1-\mu) -1/N+\mu/p & m = n. \nonumber
\end{cases}
\end{align}  

\begin{enumerate}
	\item[(\textbf{A1})] The random processes $\{w_{mn}(t)\}_{t\geq 0}$ and $\{\delta_{mn}(t)\}_{t\geq 0}$ are independent identically distributed (i.i.d.). 
	\item[(\textbf{A2})] The random variables $\{\delta_{mn}(t)\}$ have support $(0, C_\delta]$, while the weights $\{w_{mn}(t)\}$ have support $\{0\} \cup [\epsilon_w, 1]$. 
	\item[(\textbf{A3})] The online algorithm is initialized such that $\norm{(\I-\bm{11}^T/N)\Xh_0} \leq C_x$.
	\item[(\textbf{A4})] There exists $t_0$ such that for any $\mu \in (0,1)$, there exists $\varrho \in (0,1)$ such that $\norm{\prod_{s = \tau+1}^t \boldsymbol{\Theta}_s}_2 < \varrho^{t-\tau}$ for all $t - \tau \geq t_0$.
		\item[(\textbf{A5})] For each $t$, the non-zero weights $\{w_{mn}(t)\}_{m,n}$ are i.i.d. with $\bar{w}:=\Ex{w_{mn}(t)}$. 
\end{enumerate}

The i.i.d. assumption in (\textbf{A1}) is standard in the analysis of most stochastic approximation algorithms. For the applications at hand, the support of $\delta_{mn}(t)$ and $w_{mn}(t)$ is naturally finite. It is required from (\textbf{A2}) that the non-zero weights be bounded away from zero. Such a condition is required to ensure the numerical stability of the Laplacian system of equations that must be solved at every iteration [cf. \eqref{smdseq}, \eqref{smdsg}]. Specifically, it is shown in Appendix that (\textbf{A2}) implies the following result
\begin{lem}\label{lemeps}
Under (\textbf{A2}), it holds that $\norm{\L_t^\dagger}_2 \leq \epsilon_{\L} := (N-1)^2/2\epsilon_w $ for all $t\geq 1$. 
\end{lem}
The proof of Lemma \ref{lemeps} is provided in Appendix \ref{lemepsproof}. The initial configuration can always be normalized to satisfy the bound in (\textbf{A3}). Assumption (\textbf{A4}) restricts the extent to which the graphs $\G_t$ can stay disconnected over time. To obtain intuition on  (\textbf{A4}), observe first the largest eigenvalue of $\boldsymbol{\Theta}_t$ is $1-\mu$ if $\G_t$ has a single connected component and one otherwise. Consequently, if all $\{\G_s\}_{s = \tau+1}^t$ are connected, (\textbf{A4}) holds with $\varrho = 1-\mu$. Conversely, it holds that $\norm{\prod_{s = \tau+1}^t \boldsymbol{\Theta}_s}_2 = 1$ if and only if (a) each $\{\G_s\}_{s = \tau+1}^t$ has more than one components, and (b) the components do not change over time, i.e., $j_m(t) = j_n(t)$ for all $m$, $n$, and $t$. Intuitively, (\textbf{A4}) allows $\{\G_s\}$ to have multiple connected components at each $s \geq 1$, as long as the nodes belonging to these components keep changing over time. 

Finally, (\textbf{A5}) is perhaps the most restrictive, and may not always be easy to satisfy. For instance, the weights are not identically distributed in the context of dynamic network localization (cf. Sec. \ref{dnl}), since non-zero weights are often assigned to neighboring nodes only. Likewise, weights selected via Sammon mapping also result in non-identically distributed weights. The assumption however greatly simplifies the proof of convergence for the averaged algorithm. Having stated the assumptions, the averaging analysis is presented in the subsequent subsection. 

\subsubsection{Hovering Theorem}
The proposed stochastic SMACOF algorithm will be related to an averaged algorithm with updates,
\begin{align}\label{avgd}
\Xt_{t+1} = (1-\mu\upsilon)\Xt_t + \mu \Ba(\Xt_t)\Xt_t 
\end{align}
where the time-invariant function $\Ba(\X) := \Ex{\L_t^{\dagger}\Bt(\X)}$ and $\upsilon = \frac{N(p-1)}{p(N-1)}$. Assuming that both algorithms start from the same initialization, i.e.,  $\Xh_0 = \Xt_0$, the following proposition states the main result of this section.

\begin{prop}\label{hovstat}
Under (\textbf{A1})-(\textbf{A5}), and for $\mu<1$, it holds for the updates generated by \eqref{smdsg} and \eqref{avgd}, that
\begin{align} \label{hoveq}
\max_{1\leq t\leq 1/\mu} \norm{\Xh_t-\Xt_t} \leq c(\mu)
\end{align}
\vspace{-.5mm}
where the random variable $c(\mu) \rightarrow 0$ almost surely as $\mu \rightarrow 0$ with probability 1.
\end{prop}

Intuitively, Proposition \ref{hovstat} states that the trajectory of the proposed stochastic algorithm in \eqref{smdsg} stays close to that of the averaged algorithm in \eqref{avgd}. Further, the stochastic ''oscillations'' of \eqref{smdsg} are small if $\mu$ is also small. However, choosing too small a value of $\mu$, which is also the step-size in \eqref{avgd}, will generally result in a slower convergence rate for any such iterative algorithm. The parameter $\mu$ may therefore be seen as controlling the trade-off between the convergence rate and asymptotic accuracy. Further characterization of this trade-off is pursued via numerical tests in Sec. \ref{sim}.   

Alternatively, consider the case when $T$ updates of \eqref{smdsg} are performed with $\mu = 1/T$. For this case, the bound in \eqref{hoveq} becomes
\begin{align}
\max_{1\leq t\leq T} \norm{\Xh_t-\Xt_t} \leq c(1/T) \label{hoveqT}
\end{align}
where $c(1/T) \rightarrow 0$ almost surely as $T \rightarrow \infty$. In other words, the stochastic oscillations can be made arbitrarily small if sufficient number of updates can be performed. It is remarked that such results are commonplace in the stochastic approximation literature \cite{solo1994adaptive,kushner2003stochastic,borkarstochastic}. 

Next, an outline of the proof of Proposition \ref{hovstat} is presented, while the details are deferred to Appendix \ref{lemproof}. The overall structure of the proof is similar to that in \cite{solo1994adaptive}. Significant differences exist in the details however, since workarounds are introduced in order to avoid making any assumptions on the boundedness of $\Xh_t$. It is emphasized that such a modification is generally not possible in a vast majority of problems, and is not trivial. It is however possible here due to the special structure of the update \eqref{smdsg} that depends only on the differences between pairs of rows of $\Xh_t$; see \eqref{mrow}.

\begin{IEEEproof}[Proof of Proposition \ref{hovstat}]
The difference between the iterates generated by \eqref{smdsg} and \eqref{avgd} is given by
\begin{align}
\Dl_{t+1}  := \Xh_{t+1}-\Xt_{t+1}& = \Dl_t -\mu \left(\L_t^{\dagger}\L_t\Xh_t - \upsilon\Xt_t \right) \nonumber\\
&\hspace{-1cm}+ \mu \left(\L_t^\dagger\B^\epsilon_t(\Xh_t)\Xh_t - \Ba(\Xt_t)\Xt_t \right) \label{dlt2}
\end{align}
Assuming that both the algorithms start from the same initialization, i.e.,  $\Xh_0 = \Xt_0$, it follows that
\begin{align}
\Dl_{t+1} &=  -\sum_{\tau=0}^t\mu \left(\L_\tau^{\dagger}\L_\tau\Xh_\tau - \upsilon\Xt_\tau \right) \nonumber\\
&+ \mu \sum_{\tau=0}^t\left(\L_\tau^\dagger\B^\epsilon_\tau(\Xh_\tau)\Xh_\tau - \Ba(\Xt_\tau)\Xt_\tau \right) \nonumber\\
&= -\mu \upsilon\sum_{\tau = 1}^t\Dl_{\tau} + \mu \left(\K^{1}_t + \K^2_t + \K^3_t\right)
\end{align}
where for all $t \geq 0$, 
\begin{subequations}\label{kt}
\begin{align}
\K^1_t &= \sum_{\tau=0}^{t} \left(\L^{\dagger}_\tau\B^\epsilon_\tau(\Xh_\tau)\Xh_\tau - \Ex{\L^{\dagger}_\tau\B^\epsilon_\tau(\Xh_\tau)}\Xh_\tau\right) \label{kt1}\\
\K^2_t &= -\sum_{\tau=0}^{t}\left(\L_\tau^{\dagger}\L_\tau - \upsilon\I\right)\Xh_\tau\label{kt2} \\
\K^3_t &= \sum_{\tau=1}^{t}\left( \Ba(\Xh_\tau)\Xh_\tau - \Ba(\Xt_\tau)\Xt_\tau\right)\label{kt3}.
\end{align}
\end{subequations}

The following intermediate lemma develops bounds on the three terms in \eqref{kt}, and constitutes the key step in the proof. 

\begin{lem}\label{ktlem}
The following bounds hold for $i = 1,2$
\begin{align}
\norm{\K_t^i} &\leq d^i_t + C_i \mu \sum_{\tau = 1}^{t} \pi^i_t \label{ktib}\\
\norm{\K_t^3} &\leq C_3 \sum_{\tau = 1}^{t} \norm{\Dl_\tau} 
\end{align}
where the constants $C_1$, $C_2$, and $C_3$ are independent of $t$, and the constants $d^1_t$, $\pi^1_t$, $d^2_t$, and $\pi^2_t$ are such that
\begin{subequations}\label{dpit}
\begin{align}
\frac{d^i_t}{t} &\rightarrow 0 & \frac{\pi^i_t}{t} &\rightarrow 0
\end{align}
for $i = 1, 2$, almost surely as $t \rightarrow \infty$. 
\end{subequations}
\end{lem}

The proof of Lemma \ref{ktlem} is provided in Appendix \ref{lemproof}. The norm of $\Dl_{t+1}$ can therefore be bounded by applying triangle inequality on \eqref{dlt2} as follows. 
\begin{align}\label{dltfinal}
\norm{\Dl_{t+1}} \leq \mu (C_3+1)\sum_{\tau = 1}^t \norm{\Dl_\tau} + f(\mu)
\end{align}
where we have used the fact that $\norm{\J}_2 = 1$ and
\begin{align}
f(\mu) := \max_{0\leq t\leq 1/\mu}\mu (d^1_t + d^2_t) + \mu^2\sum_{\tau=1}^t C_1\pi^1_t + C_2\pi^2_t.
\end{align}
It is further shown in Appendix \ref{lemproof} that $f(\mu) \rightarrow 0$ almost surely as $\mu \rightarrow 0$. Proposition \ref{hovstat} then follows from the application of the discrete Bellman-Gronwall Lemma \cite{solo1994adaptive} on \eqref{dltfinal}, which yields
\begin{align}
\norm{\Dl_t} & \leq f(\mu) (1+\mu(C_3+1))^{t} = f(\mu)e^{t\log(1+\mu(C_3+1))} \nonumber\\
&\leq f(\mu)e^{\mu t(C_3+1)} \leq f(\mu)e^{C_3+1} := c(\mu) 
\end{align}
\end{IEEEproof}

\subsubsection{Convergence of the Averaged Algorithm}
Having established that the trajectory of the stochastic algorithm hovers around that of the averaged algorithm, we complete the proof by establishing that the averaged algorithm converges to a local minimum of \eqref{staticmds}. The challenge here is that the updates in \eqref{avgd} do not resemble those in other classical algorithms such as SMACOF or gradient descent. For notational brevity, let $\bar{\delta}_{mn} := \Ex{\delta_{mn}(t)}/\sqrt{\norm{\x_m-\x_n}^2+\epsilon_x}$ and $\J = \I - \b{11}^T/N$, and note the following result.

\begin{lem}\label{avglem}
Under (\textbf{A1})-(\textbf{A5}), it holds that 
\begin{align}
\left[\Ba(\X)\right]_{mn}= \frac{\upsilon}{N}\begin{cases}  - \bar{\delta}_{mn} & m\neq n \\
\sum\limits_{k\neq m} \bar{\delta}_{mk} &m = n.
\end{cases}\nonumber
\end{align}
\end{lem}
The proof of Lemma \ref{avglem} is provided in the Appendix \ref{avglemProof}. For the rest of the section, we will assume that $\epsilon_x \ll 1$ and thus negligible. Therefore from Lemma \ref{avglem}, we have that 
\begin{align}
\bar{\sigma}(\X) &= \sum_{m<n} \Ex{w_{mn}(t)\delta_{mn}(t)^2} + \tr{\X^T\bar{\L}\X} \nonumber\\
& ~~ - 2\text{tr}(\X^T\bar{\B}(\X)\X) 
\end{align}
where, $\bar{\L} := \Ex{\L_t} = \bar{w}\upsilon p \J$ and 
\begin{align}
\bar{\B}(\X) &:= \frac{\bar{w}\upsilon p}{N}\begin{cases}  - \bar{\delta}_{mn} & m\neq n \\
\sum\limits_{k\neq m} \bar{\delta}_{mk} &m = n. 
\end{cases}
\end{align}
The main result of this subsection is stated as the following proposition. 

\begin{prop}\label{avgconv}
The mean-stress values $\bar{\sigma}(\Xt_t)$ decrease monotonically with $t$ and converge to a stationary point of \eqref{ss}.
\end{prop}
\begin{IEEEproof}
Without loss of generality, let $\sum_{m<n} \Ex{w_{mn}(t)\delta_{mn}(t)^2} = 1$, and define $\eta^2(\X):=\frac{1}{\mu\upsilon}\tr{\X^T\bar{\L}\X}$ and $\rho(\X):=\frac{1}{2}(1/\mu\upsilon-1)\tr{\X^T\bar{\L}\X} + \tr{\X^T\bar{\B}(\X)\X}$, and observe that $\bar{\sigma}(\X) = 1 + \eta^2(\X) -2\rho(\X)$. Similarly, define the mapping $\bG(\X) := (1-\mu\upsilon)\X + \frac{\mu}{\bar{w}p}\bar{\B}(\X)\X$, so that the updates in \eqref{avgd} become $\Xt_{t+1}=\bG(\Xt_t)$. 

Given any two embeddings $\X$ and $\Y$, the following bounds hold from the Cauchy-Schwarz inequality:
\begin{align} \label{majo:a}
-\text{tr}(\X^T\bar{\L}\X) &\leq -\text{tr}((2\X-\Y)^T\bar{\L}\Y)  \\
-\text{tr}(\X^T\bar{\B}(\X)\X) &\leq -\text{tr}(\X^T\bar{\B}(\Y)\Y)  
\end{align}
which allows us to conclude that 
\begin{align}
\rho(\X) &\geq \frac{1}{\mu\upsilon}\tr{\X^T\bar{\L}\bG(\Y)} -\frac{1-\mu\upsilon}{2\mu\upsilon}\tr{\Y^T\bar{\L}\Y} \\
\Rightarrow ~\bar{\sigma}(\X) &\leq 1 + \eta^2(\X) + \frac{1-\mu\upsilon}{\mu\upsilon}\tr{\Y^T\bar{\L}\Y} \nonumber\\ &- \frac{2}{\mu\upsilon}\tr{\X^T\bar{\L}\bG(\Y)} \nonumber\\
&\hspace{-1.5cm}=1+(1-\mu\upsilon)\eta^2(\Y) - \eta^2(\bG(\Y)) + \eta^2(\X-\bG(\Y))  \label{rhox}
\end{align}
where equalities holds for $\X = \Y$. Denote the right-hand side of \eqref{rhox} by $\omega_{\Y}(\X)$, and observe that $\omega_{\Y}(\X) \geq \omega_{\Y}(\bG(\Y))$ for all $\X$. This yields the main inequality that $\bar{\sigma}(\Y) = \omega_{\Y}(\Y) \geq  \omega_{\Y}(\bG(\Y)) \geq \bar{\sigma}(\bG(\Y))$. In other words, we have that $\bar{\sigma}(\Xt_{t}) \geq \bar{\sigma}(\Xt_{t+1})$, so that the non-negative sequence $\sigma_t := \bar{\sigma}(\Xt_t)$ is non-increasing and therefore convergent to a limit, say $\bar{\sigma}_{\infty}$. By squeeze theorem for limits \cite{rudin1964principles}, it also holds that $\omega_{\Xt_t}(\Xt_{t+1}) \rightarrow \bar{\sigma}_{\infty}$, yielding the following limits
\begin{align}
\lim_{t\rightarrow \infty}\eta^2(\X_t) &= (1-\bar{\sigma}_{\infty})/\mu \upsilon \\
\lim_{t\rightarrow \infty}\rho(\X_t) &= (1-\bar{\sigma}_{\infty})(1-\mu \upsilon)/2 \\
\lim_{t\rightarrow \infty}\eta^2(\Xt_t-\Xt_{t+1}) &= 0 \label{errlim}
\end{align}
Since the matrices $\{\Xt_t\}_{t\geq 0}$ are origin centered, the result in \eqref{errlim} can equivalently be written as $\norm{\Xt_t - \Xt_{t+1}} \rightarrow 0$ as $t \rightarrow \infty$. Denoting the limit point of $\Xt_t$ by $\Xt_{\infty}$, it can be seen that $\nabla \bar{\sigma}(\Xt_{\infty}) = 0$. 
\end{IEEEproof}

\section{Implementation Aspects}\label{imp}
\subsection{Multi-agent network localization}\label{dnl}
Multidimensional scaling has been widely used for localization, where inter-node distances are often obtained from time-of-arrival or received signal strength measurements \cite{hero-dwmds,patwari2005locating,wymeersch2009cooperative,simonetto2014distributed}. Wireless network localization is challenging because the pairwise distance measurements are noisy, time-varying due to mobility, fading, and synchronization errors, and often partially missing, due to the limited range of the sensors. Further, the limited battery life and resource constraints at the nodes impose restrictions on the communication and computational load that the network can tolerate\cite{kumar2016cooperative,simonetto2014distributed}. 

Towards addressing these limitations, the stochastic SMACOF algorithm for network localization works by judiciously choosing $\{w_{mn}(t)\}$ to limit the communication and computational cost at each update. The idea is to partition the network into several non-overlapping clusters (or components), chosen randomly at each time $t$. The coordinates within a cluster are updated as in \eqref{smdsg}. Only neighboring nodes are included within each cluster, thus eliminating the need for multihop communication between far off nodes. Finally, since the updates at different components are independent of each other, the localization algorithm is run asynchronously as follows.
\begin{enumerate}
	\item[S1. ] At a given time $t$, a node $j$ randomly declares itself as a cluster head, and solicits cluster members from among its neighbors $n \in \mathcal{N}_j$. Available neighbors respond with their current location estimates $\hat{\x}_n(t)$, resulting in a star shaped cluster $\mathcal{C}_t^j$. Once locked as cluster members, these nodes respond only to the messages from node $j$.	
	\item[S2. ] The cluster head performs distance measurements between itself and all its neighbors and collects $\delta_{jn}(t)$ for all $n \in \mathcal{C}_t^j\setminus\{j\}$. 
	\item[S3. ] The cluster head performs the update in \eqref{smdsg} with appropriately chosen weights $\{w_{jn}(t)\}$, and broadcasts the new location estimates to each node in $\mathcal{C}^j_t\setminus\{j\}$
	\item[S4. ] Nodes in $\mathcal{C}^j_t\setminus\{j\}$, upon receiving the new location estimates (or upon timeout or error events), release their locks and become available.
\end{enumerate}

As originally intended, the proposed algorithm can also be applied to mobile networks. The algorithm is expected to perform well as long as the node velocities are not too high. The asynchronous nature of the algorithm allows for delayed updates at nodes, balanced battery usage within the network, and communication errors. In general, it is also possible to apply multiple updates of the form in  \eqref{smdsg} per time instant, without incurring any extra communication cost. 

Nodes may declare themselves as cluster heads using a random backoff-based contention mechanism such as CSMA, and solicit neighbors by simply sending an RTS packet. An update at a cluster thus takes up at most two message exchanges. More complicated protocols that ensure recovery from collisions, and robustness or errors can also be used\cite{demirkol2006mac}. The online algorithm is flexible, and allows clusters of any shape or size, depending on the communication and computational resources available within the network. The non-zero weights, corresponding to available distance measurements, can be chosen according to the estimated noise variance \cite{hero-dwmds,patwari2005locating}, following Sammon mapping \cite{mds-book,torgerson1965multidimensional}, or simply as unity.


It is remarked that the node coordinates obtained from (S1)-(S4) are relative and centered at the origin. In applications where node coordinates are required with respect to a set of GPS-enabled anchor nodes, appropriate rotation and translation operations must be applied at each node. Since the anchor nodes are generally not power constrained, it is possible for them to determine these transformations \cite{hero-dwmds}, and convey the result to all other nodes. As shown later in Sec. \ref{sim}, it is generally sufficient to calculate the transformations periodically every few time slots. 

Finally, similar to the SMACOF algorithm, the stochastic SMACOF algorithm is sensitive to initialization. A random initialization may result in the algorithm getting trapped in a ``poor'' local minimum. In practice, superior location estimation performance is obtained if the initialization is at least roughly correct. Simple low-complexity localization algorithms can be used for initialization. For instance, nodes can roughly triangulate themselves using noisy distance estimates from the anchor nodes \cite{savarese2001location}. 

\subsection{Large network visualization}\label{viz}
It is possible to visualize $N$ objects in a 2 or 3 dimensional euclidean space by applying MDS to the pairwise dissimilarities $\{\delta_{mn}\}$. The SMACOF algorithm is however ill-suited for large-scale visualization since it requires at least $\O(N^2)$ operations per iteration. Further, even processing the full measurements $\{\delta_{mn}\}$ simultaneously may not be feasible for datasets with more than a hundred thousand objects. 

Visualization via stochastic embedding can be achieved by partitioning the objects into several subsets of reasonable sizes,  and performing the updates in \eqref{smdsg}. The following steps are performed for each $t \geq 1$. 
\begin{enumerate}
	\item Partition the $N$ objects into random, mutually exclusive subsets $\cjt$ with $p$ nodes per subset.
	\item For each subset, randomly choose a small fraction $f_t$ of pairs and measure (calculate or fetch from memory) distances $\delta_{mn}$ for the chosen pairs. { Let $\mathcal{F}_t^j$ denote the set of chosen pairs for each cluster $j$ and time $t$.}
	\item Apply the update in \eqref{smdsg} for each subset $\mathcal{C}_t^j$. 
\end{enumerate}
Compared to the localization algorithm, in this case all pairwise distances are available a priori and without noise, but cannot be read or processed simultaneously. The aforementioned steps result in making $\{w_{mn}(t)\}$ sparse and thus reducing the per-iteration complexity. 
{  Algorithm \ref{smacof:visualization} summarizes the implementation of stochastic SMACOF for large network visualization.}

\begin{algorithm}
 \caption{Stocahstic SMACOF for Large Network Visualization}\label{smacof:visualization}
\begin{algorithmic}[1]
{ 
\STATE Initialize $\X_{0}$ and set $\mu$ to some value in $(0,1)$
 \FOR{$t=1,2,\ldots$}
\STATE Partition the set $\N$ into $C$ disjoint subsets $\{\Cs_t^j\}_{j=1}^C$  
\FOR{$j=1,\ldots ,C$}
\STATE Measure or fetch from memory pairwise distances $\{\delta_{mn}(t)\}$, for a subset of object pairs $(m,n) \in \mathcal{F}_t^j$.
\STATE Set weights $w_{mn}(t)=1$ for all $(m,n) \in \mathcal{F}_j^t$. 
\STATE Perform the update in \eqref{smdsg} for each subset $\mathcal{C}_t^j$.
\ENDFOR
\ENDFOR
}
 \end{algorithmic}
 \end{algorithm}

Again, as envisioned earlier, the algorithm is also applicable to visualization of dynamic networks. The idea here is to create an animation consisting of embeddings that vary over time. By specifying a small enough value for $\mu$ in \eqref{smdsg}, it is possible to force the embeddings to change slowly over time, thus preserving the user's \emph{mental map} \cite{xu2013regularized}. Unlike existing algorithms however, the proposed algorithm can allow visualization of very large datasets.

\subsection{Algorithm complexity}
{  Unlike the SMACOF algorithm, whose per-iteration complexity is $\O(N^2)$, the stochastic SMACOF algorithm processes the data in small batches and can therefore be implemented at near-linear complexity.} This is because if $\G_t$ consists of multiple components of size $p$ each, the updates in \eqref{smdsg} decouple and can even be carried out in parallel. Further, the weights for each cluster are chosen to be sparse, i.e., the $p \times p$ matrix $\L_t^{(j)}$ has at most $q \ll p^2$ non-zero elements. The problem of solving a sparse Laplacian system of equations  has been well studied, and state-of-the-art solvers return a solution in time $\O(q \log p)$ for each component. Thus, using $N/p$ sparse matrices $\{\L_t^j\}$ results in an overall complexity of $\O\left(\frac{Nq}{p}\log p\right)$. {  As we will show next, the appropriate choice of the batch size $p$ results in a near-linear complexity. The complexity results obtained in this section are summarized in Table \ref{complex}. }

Note that a sublinear per-iteration complexity of $\O(q\log(p))$ is also achievable by updating only one component per iteration. Such an implementation would however require proportionally large number of iterations. Alternatively, the per-iteration complexity of the algorithm can be calibrated using the total number of dissimilarity measurements processed per-iteration, given by $f(N) = q(N/p)$. To this end, we provide approximate rules for choosing $p$ and $q$ so as to minimize the per-iteration complexity, given the total number of non-zero weights $f(N)$.


First, assume that each $\L_t^j$ is sparse, i.e., $q \ll p^2$, so that $f(N)/N = q/p  \ll p$. In this case, since the per-iteration complexity is given by $\O(f(N)\log(p))$, the value of $\log(p)$ should be as small as possible. It can be seen that the choice 
\begin{align}
p &\sim \O\left(\left(\frac{f(N)}{N}\right)^\beta\right) & q &\sim \O\left(\left(\frac{f(N)}{N}\right)^{\beta+1}\right)
\end{align}
for some $\beta \gg 1$ results in the complexity $\O(f(N)\log(f(N)/N))$, while ensuring that $\L_t^j$ is still sparse with $q \sim \O(p^{1+1/\beta})$. Note that it is not necessary for $\beta$ to be very large, as long as the sparse Laplacian solvers can still be used. On the other hand, when $\L_t^j$ is dense so that $q \sim \O(p^2)$, the per-iteration complexity is given by $\O(Nq) = \O(f(N)p)$. In this case, it holds that $f(N)/N = q/p \leq p$, so that one must choose $p \sim \O(f(N)/N)$ and $q \sim \O(f(N)^2/N^2)$. Consequently, the optimal iteration complexity for this case becomes $O(f(N)^2/N)$.

Table \ref{complex} shows a few example choices of $f(N)$ and the corresponding per-iteration complexity values. It can be observed that when $f(N)$ is almost linear in $N$, so is the per-iteration complexity, regardless of the sparsity of $\L^j_t$. On the other hand, using a sparse $\L_t^j$ becomes important when $f(N)$ is large. 

\begin{table}[htbp]
\centering
	\begin{tabular}{|c | c | c|}
	\hline
  Non-zero weights $f(N)$ & sparse $\L_t^j$ & dense $\L_t^j$ \\
	\hline
	$\O(N^{1+\omega})$, $0 < \omega \ll 1$ & $\O(N^{1+\omega}\log(N))$ & $\O(N^{1+2\omega})$ \\
	$\O(N\log N)$ & $\O(N \log N\log\log N)$ &  $\O(N \log^2 N)$ \\
	$\O(N^{3/2})$ & $\O(N^{3/2}\log N)$ & $\O(N^2)$ \\
	\hline 
\end{tabular}
	\caption{ Algorithm complexity for different choices of $f(N)$}
	\label{complex}
\end{table}


 \vspace{-7mm}

\section{Simulation results} \label{sim}
This section provides simulation results evaluating the performance of the proposed algorithm. The general properties of the stochastic SMACOF algorithm are first characterized using numerical tests. Next, simulation results are provided for the online localization algorithm, evaluating its performance in various mobile network scenarios. Finally, applicability to large-scale visualization is demonstrated by running the algorithm on two different datasets. {  Before proceeding, it is remarked that the proposed stochastic SMACOF is better suited to applications where the size of the dataset is large, preferably $N>50$. Indeed, if the problem at hand is small (say $N<20$), conventional SMACOF would likely be faster, since the proposed algorithm generally requires more iterations to converge. The computational advantage arising from processing only a few distance measurements per time instant becomes significant only when $N$ is sufficiently large.}

\vspace{-3mm}
\subsection{Algorithm Behavior}
This section provides several numerical tests that allow us to study various properties of the stochastic SMACOF algorithm. Towards this end, consider a network with 100 nodes, distributed uniformly over a $10 \times 10$ planar area. The measured distances between nodes $m$ and $n$ are given by $\delta_{mn}(t) = \norm{\x_m-\x_n} + v_{mn}(t)$, where $v_{mn}(t) \sim \mathcal{N}(0,0.01)$. Negative distance measurements were discarded by setting the corresponding $w_{mn}(t) = 0$. The algorithm is run for different values of $\mu$, with $p = 25$ and about $35\%$ density of non-zeros\footnote{Non-zero locations are generated randomly, and the number of non-zeros vary between different instantiations.}. All non-zero weights are chosen to be unity.


\begin{figure}[t]
  \centering
  \includegraphics[scale=1]{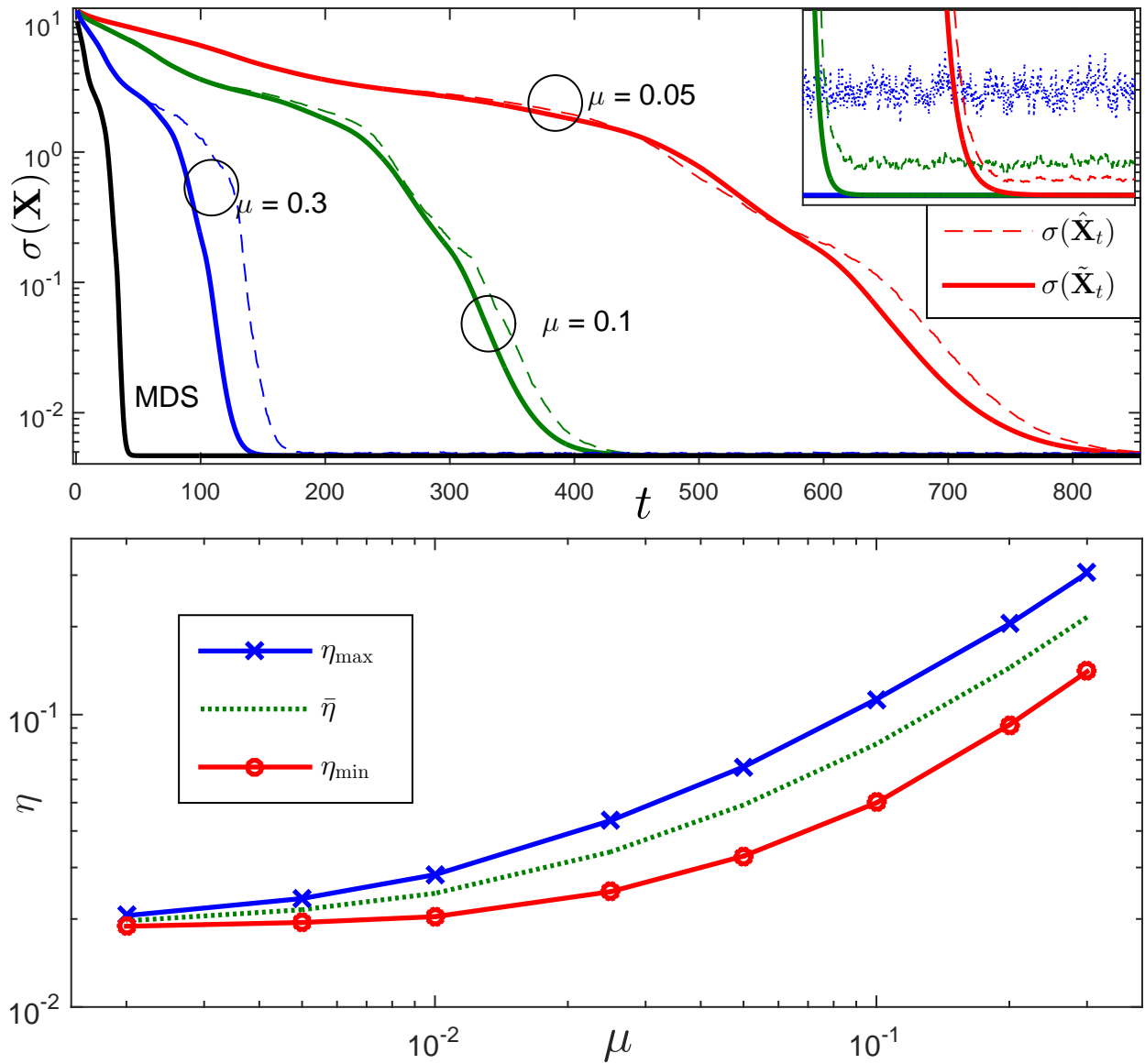}
  \vspace{-.5cm}
  \caption{(Top) Performance of the stochastic SMACOF algorithm, the averaged algorithm, and the SMACOF algorithm; (Bottom) Steady state fluctuations in the stress.  } 
  \label{avg_hov_combined}
\end{figure}

\subsubsection{Transient performance}
 Fig  \ref{avg_hov_combined} (Top) shows the sequence of normalized stress values obtained from an example run of the algorithm [cf. \eqref{smdsg}]. For comparison, the stress values obtained from running the averaged algorithm (cf. \eqref{avgd}) and the SMACOF algorithm for weighted MDS (cf. \eqref{smacof}) are also plotted. All algorithms are intialized with the same randomly chosen configuration. The MDS algorithm runs with all-one weights, while the updates for the averaged algorithm are obtained via empirical averaging. 

As expected, the convergence speed of the algorithm varies monotonically with $\mu$. Consistent with Proposition \ref{hovstat}, the trajectory of the proposed algorithm follows that of the averaged algorithm. As expected, the steady-state stress value achieved by the averaged algorithm is very close to that of SMACOF. Further, as shown in the inset, the proposed algorithm hovers above the averaged algorithm, with steady-state deviation decreasing with $\mu$. 

It is remarked that the SGD algorithm, with updates specified in \eqref{sgd}, tended to diverge in the presence of noisy distance measurements, different weight choices, and poor initializations. For instance, when using Sammon mapping, i.e., $w_{mn} = 1/\delta_{mn}$, the noisy measurement model specified earlier, and $\mu = 0.05$, the SGD algorithm converged for only 19 out of 100 test runs. In contrast, no divergent behavior was ever observed for the proposed algorithm even with measurement noise $v_{ij}\sim\mathcal{N}(0,10)$.

%

\subsubsection{Steady state performance}
The algorithm is allowed to run for 5000 time instants with different values of $\mu$, and the minimum, mean, and maximum steady-state stress values are evaluated. We set $\mathcal{T}_{\text{ss}} = [4801, \ldots, 5000]$ and evaluate 
\begin{align}
\eta_{\min} &= \min_{t\in \mathcal{T}_{\text{ss}}} \sigma(\Xh_t) & \bar{\eta} &=  \sum_{t\in \mathcal{T}_{\text{ss}}}\frac{\sigma(\Xh_t)}{\abs{\mathcal{T}_{\text{ss}}}} &
\eta_{\max} &= \max_{t\in \mathcal{T}_{\text{ss}}} \sigma(\Xh_t). \nonumber
\end{align}
Starting with the same initialization, the entire experiment is repeated for 100 Monte Carlo iterations. Fig. \ref{avg_hov_combined} (Bottom) shows the minimum, mean, and maximum steady state errors plotted against $\mu$. As expected, the stress values converge to a small non-zero value that decreases with $\mu$.

\vspace{-3mm}

%
%
%
%
%
%
\subsection{Dynamic Network Localization}
The localization performance of the proposed algorithm is studied on a mobile network.  Video\footnote{https://www.youtube.com/watch?v=-MQFR3yiv7U } shows an example run of the algorithm on a mobile network with $N = 8$ and $\mu = 0.3$. The performance of the algorithm is further analyzed by carrying out simulations over networks with different sizes and node velocities. For a mobile network with $N$ nodes, nodes are deployed randomly with an average density of one node per unit area. Nodes can measure distances and communicate within a radius of $\sqrt{N}/2$. For all values of $N$, five nodes are randomly chosen to be anchors. The node velocities are initialized randomly and updated according to the following model $\mathbf{v}_{mn}(t+1) = \alpha\mathbf{v}_{mn}(t) + \sqrt{1-\alpha^2}\mathbf{n}_v(t)$, where $\mathbf{v}_{mn}(0), \mathbf{n}_v(t) \sim \mathcal{N}(0,\sigma_v^2\I)$. The mobility parameter $\sigma_v$ is directly proportional to the average speed of the nodes, and influences the tracking performance of the algorithms used.

The performance of the proposed algorithm is compared with the weighted MDS solution obtained by running the SMACOF algorithm till convergence. The non-zero weights, corresponding to node pairs within the communication radius of each other, are all set to one. Note however that a direct comparison between the SMACOF  solution and the proposed algorithm is unfair, since SMACOF is too complex to be directly implemented in a mobile network. Even among cooperative localization techniques that focus on efficient implementation (see e.g. \cite{hero-dwmds,patwari2005locating,wymeersch2009cooperative,simonetto2014distributed}), localization requires several iterations per time instant. In contrast, the proposed algorithm is asynchronous, and incurs linear or sublinear complexity, but is inaccurate for the first few time instants. 

In order to perform a fair comparison between algorithms, the following modifications are adopted. First, a time-slotted version of the stochastic SMACOF algorithm is considered. Within each time slot, the network forms several clusters, and performs steps (S1)-(S4). In order to reduce the overhead associated with cluster formation, nodes with fewer than 5 neighbors do not form clusters. Similarly, to limit the computational complexity at each node, cluster heads respond to at most 10 nearest neighbors. With these settings, the computational and communication complexity incurred by the network at every time slot is approximately $N/5$. The computational and communication complexity of the SMACOF variants in \cite{hero-dwmds,simonetto2014distributed,dong2012cooperative} is also normalized appropriately. As a first order approximation, it is assumed that these algorithms require $\mathcal{O}(N)$ message exchanges per iteration. Equivalently, if we allow $N/5$ message exchanges per iteration, and assume that 10 iterations are required for convergence, SMACOF requires about 50 time slots for convergence. For obtaining the plots however, SMACOF is run till convergence, and the number of iterations incurred was often more than 50. Both algorithms start with an initial estimate of the node locations. Approximate node estimates can be quickly obtained using simple techniques such as those in \cite{savarese2001location}. For the purpose of simulations, the initial locations are chosen as $\hat{\x}_{m}(0) = \x_m(0) + \mathbf{v}_{m}$ where $v_{m} \sim \mathcal{N}(0,N/100)$. Warm starts are utilized at subsequent time slots by initializing SMACOF with the previously estimated node locations. 

\begin{figure}[t]
  \centering
  \includegraphics[scale=1]{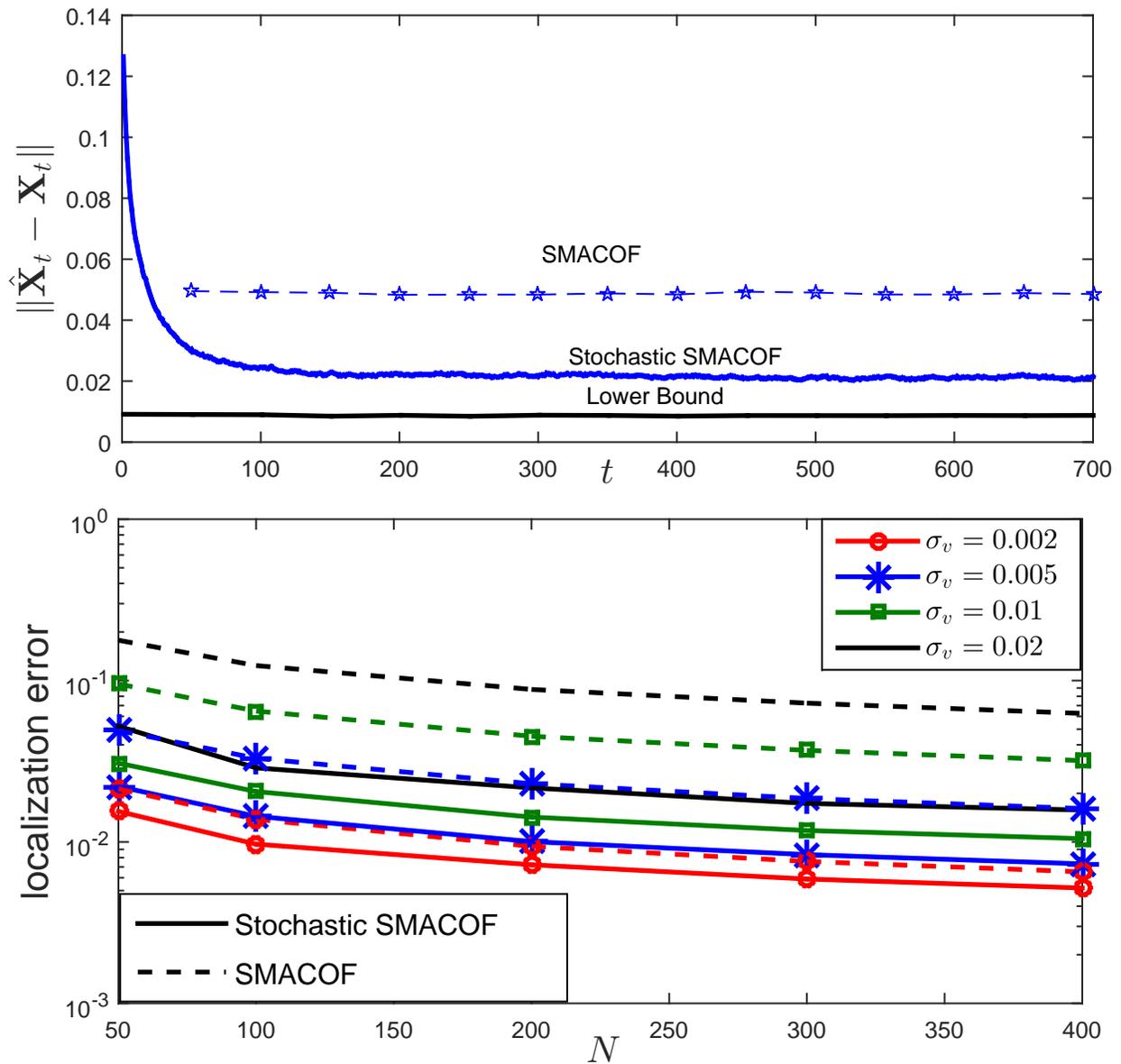}
  \vspace{-.5cm}
  \caption{(Top) Estimation error for an example run of the Stochastic SMACOF and SMACOF algorithms; (Bottom) Localization error for different network sizes and average node velocities.   } 
  \label{vssamcof_locerrvsN2_combined}
\end{figure}

Fig. \ref{vssamcof_locerrvsN2_combined}(top) shows an example run of the two algorithms with $\sigma_v = 0.01$, $N = 50$, and $\mu = 0.5$. The best possible estimation error obtained by solving the MDS problem is also shown for comparison. Observe that the proposed algorithm is inaccurate initially, and gradually approaches its steady state value. Interestingly,  the transient period required by the proposed algorithm is small, especially when compared to the 50 time slots required by the SMACOF implementation.

Next, the steady-state localization error of the two algorithms is compared. Both algorithms are run for 700 iterations, and the maximum localization error incurred in the last 200 iterations is evaluated as $e_\ell = \max_{t\in\mathcal{T}_{\text{ss}}} \frac{1}{N} \norm{\Xh_t-\X_t}$ where $\mathcal{T}_{\text{ss}} = [501, \ldots, 700]$. The entire process is repeated for 100 Monte-Carlo repetitions. For the proposed algorithm, the value of $\mu$ is tuned a priori to minimize the localization error. Fig. \ref{vssamcof_locerrvsN2_combined}(bottom) shows the steady-state localization error incurred by the online and SMACOF algorithms, plotted for different values of $N$ and $\sigma_v$. It is evident that the proposed algorithm performs significantly better than the complexity-normalized SMACOF. In particular, while the performance of the two algorithms deteriorates with increasing node mobility, the gap between their performance also increases. This is because at higher node speeds, the node locations change significantly within the 50 time slots required by SMACOF to run. Observe that for a given average node velocity, the performance of all algorithms appears to improve with $N$. However, this is simply because the average node distances increase with $N$, thereby reducing the relative average node speeds. 

\begin{figure}[t]
  \centering
  \includegraphics[scale=1.5]{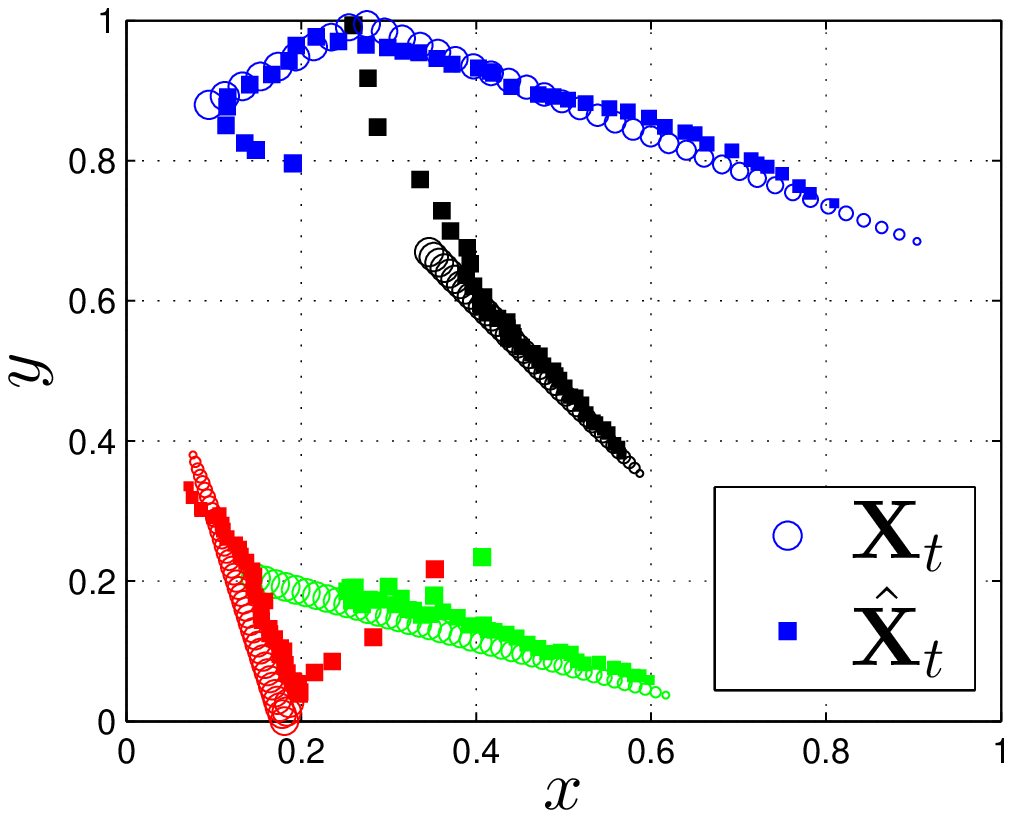}
  \caption{Example run of the dynamic network localization algorithm. Marker size decreases with time to indicate the direction of motion.}
  \label{mobilerun}
\end{figure}

{ 
  Fig. \ref{mobilerun} shows an example run of the algorithm on a mobile network with $N = 8$ and $\mu = 0.3$. The network has four static anchors placed at the four corners of the $1 \times 1$ region, that provide the necessary translation and rotation information to all other nodes. For simplicity, only one 8-node cluster is formed at each time instant by a randomly selected node. The actual and estimated node locations are shown as circles and squares respectively, with markers drawn every 10 time instants. The nodes move in the direction indicated by decreasing marker sizes. As evident from the figure, the trajectory of the estimated node locations converges to the actual trajectory within 30-40 time instants, and follows it thereafter.  }

\subsection{Large-scale Visualization}
This section demonstrates the use of the stochastic SMACOF algorithm for large-scale visualization. Given the plethora of highly sophisticated visualization algorithms a full-fledged comparison is beyond the scope of the present work. Instead, we only present the visualizations obtained from running the proposed algorithm for both static and dynamic datasets. The proposed algorithms are implemented in MATLAB  and run on an Intel Core i7 CPU. This is in contrast to the state-of-the-art visualization algorithms that require large compute clusters with hundreds of processors for similar-sized datasets \cite{baea2012visualization}.

\hspace{-1cm}\begin{figure}[htb]
  \centering
 \hspace{-1cm} \includegraphics[clip,trim={1.2cm .3cm 0cm 0 cm},scale=0.80 ]{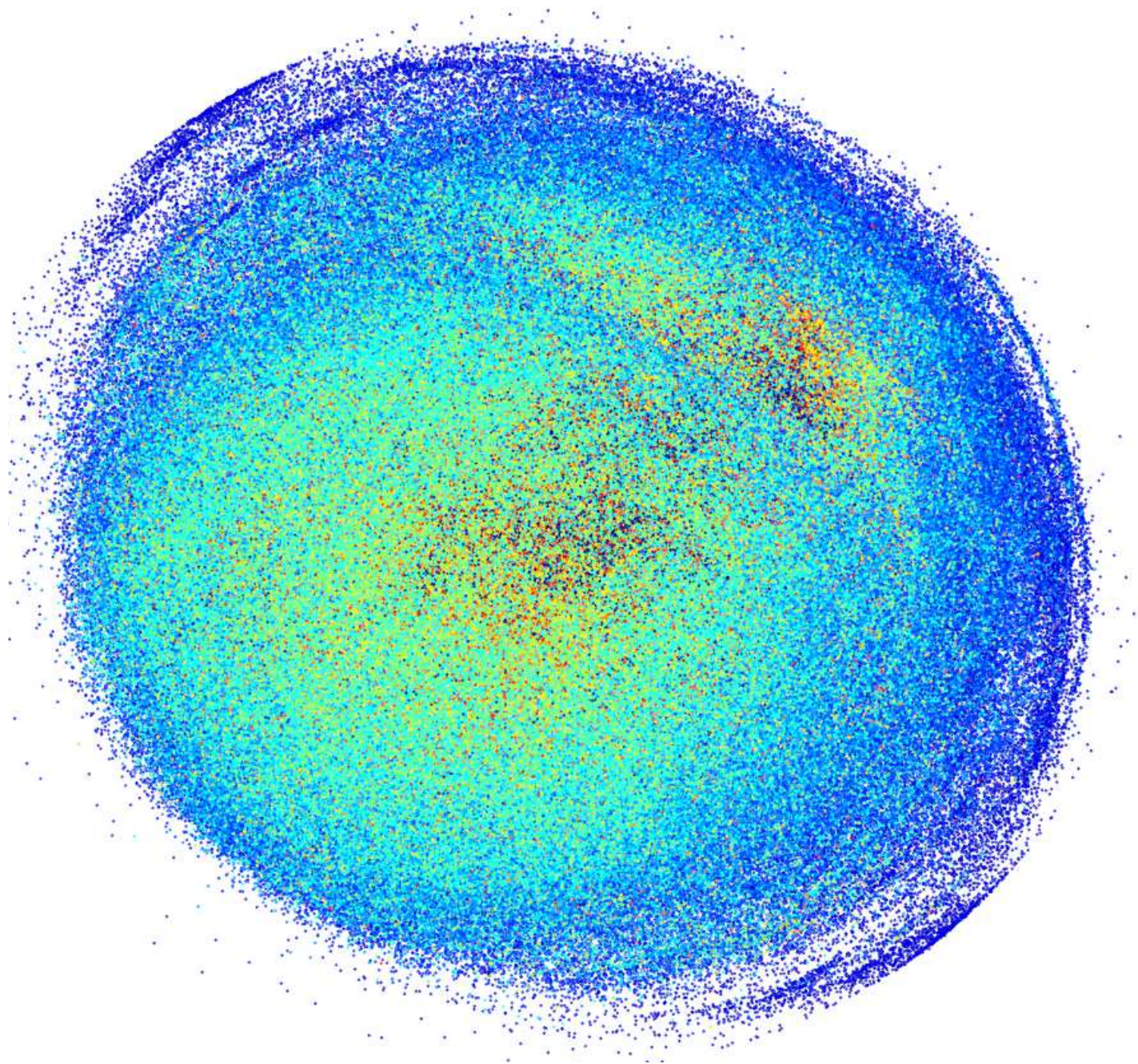}
  \caption{Visualization of  PubChem Datasets.}
  \label{pubchem}
\end{figure}

\subsubsection{PubChem Dataset} 
{ We consider a subset of 800,000 unique chemical compounds taken from the pubchem compound database\cite{kim2015pubchem ,bolton2008pubchem}. The structural information about each compound is represented by its 166 bit MACCS fingerprint. Dissimilarities between two compounds is calculated using the Tanimoto score. Dissimilarities between two compounds with binary fingerprints $\mathbf{h}$ and $\mathbf{g}$ is calculated using the Tanimoto score \cite[Ch-8]{tan2006introduction}, given by 
\begin{align}
\gamma = 1- \frac{\sum_i h_i \cap  g_i }{\sum_i h_i \cup g_i}
\end{align}
where $\cap$ and $\cup$ denote the logical AND and OR operators respectively. It is remarked that for this case, it is no longer possible to load an $N \times N$ matrix in the memory. Following the discussion in Sec. \ref{viz}, we use $p = 100$ and $q = 50$, so as to obtain linear complexity per iteration. The simulation is run for 5000 iterations, and the value of $\mu$ is reduced every 1000 iterations from 0.2 to 0.001. Figure. \ref{pubchem} shows the visualization obtained from the stochastic SMACOF algorithm. Each dot represents a compound, and is colored according to its \emph{molecular complexity}, a
measure available from the PubChem dataset. Specifically, the blue dots represent simpler (lower complexity) molecules, while green, yellow, and red colored dots represent progressively more complex molecules. It is observed that MDS yields two distinct clusters of compounds, while the lower complexity compounds are scattered towards the edges. The visualization obtained here is comparable to those obtained in \cite{choi2011browsing,baea2012visualization}.}

\subsubsection{MovieLens Dataset}
The proposed algorithm is used to perform dynamic visualization of the 27,000 movies on the MovieLens database \cite{harper2016movielens}. To this end, the time-stamp associated with each movie rating is utilized to generate a dynamic network $\G_t$ that only contains the movies released and rated till the week  number $t$. The distance between two movies is estimated from their cosine similarities.  Video shows a visualization of the evolution of the movie-space over the duration 1995-2015. Each movie is colored in accordance with its popularity, and the newly released movies start at the origin. From the video, it can be seen that the popular movies move quickly (within few weeks) towards the edge of the graph, while the less popular ones tend to remain near the center. See the video at the link \footnote{https://www.youtube.com/watch?v=iJbY3HPHAUM}.

\subsubsection{Newcomb Fraternity's Dataset}
The dynamic visualization of the Newcomb Fraternity dataset \cite{newcomb1961acquaintance} is considered. Since the dataset consists of only 16 nodes, and yields only 14 snapshots overall, computational complexity is not an issue. Nevertheless, the dynamic visualization is obtained so that it may be compared with the regularized MDS technique of \cite{xu2013regularized}.  Video\footnote{https://www.youtube.com/watch?v=G9geUI3U7Tw\&feature=youtu.be } shows the dynamic visualization obtained from running the stochastic SMACOF algorithm for 50 iterations per time slot with $\mu = 0.2$. The video is generated following the procedure similar to that in \cite{xu2013regularized}. The resulting video is quite similar to the one obtained via the graph-regularized framework of \cite{xu2013regularized}. Intuitively, the momentum term in the updates in \eqref{smdsg} plays the role of the regularization term here, and keeps the embeddings from changing too quickly.

\section{Conclusion}\label{conclusion}
The multidimensional scaling (MDS) problem is considered within a stochastic setting, and a novel stochastic scaling by majorizing a complicated function (SMACOF) is proposed. The proposed algorithm is highly scalable, and is applicable to visualization and localization problems of very large sizes. Asymptotic analysis of the stochastic SMACOF algorithm shows that it stays close to the trajectory of an averaged algorithm, which itself converges to a stationary point of the stochastic stress minimization problem. Implementation details, as well as the computational complexity analysis of the proposed algorithms are also provided. The performance of the proposed algorithm is discussed for large-scale localization and visualization examples. The efficacy of the proposed algorithm is demonstrated for localization of mobile networks, and visualization of both, static and dynamic networks.

\vspace{-3mm}
\appendices
\section{Lower bound on the algebraic connectivity} \label{lemepsproof}
In order to obtain intuition on (\textbf{A3}), consider the undirected graph $\G_t$ whose edges have weights $\{w_{mn}(t)\}$, and recall that $\L_t$ is the graph Laplacian of $\G_t$. The eigenvalues of $\L_t$ constitute the spectrum of the graph $\G_t$\cite{mohar1992laplace}. If $\G_t$ is connected, the smallest eigenvalue of $\L_t$ is zero, while the second-smallest eigenvalue $a(\G_t) = 1/\norm{\L_t^\dagger}_2$ is always non-zero and is referred to as the algebraic connectivity of $\G_t$. As the name suggests, $a(\G)$ captures the overall connectivity of the graph. On the other hand, if $\G_t$ has $K \geq 2$ connected components $\{\G_t^k\}_{k=1}^K$, the $K$ smallest eigenvalues of $\L_t$ are zero, so the smallest non-zero eigenvalue is simply $a(\G_t) = \min_k a(\G_t^k)$. Next, we establish a lower bound on the algebraic connectivity of the weighted graph $\G_t$. 

\begin{IEEEproof}[Proof of Lemma \ref{lemeps}]
If $\G_t$ is connected, the second smallest eigenvalue is given by 
\begin{align}\label{algdef}
a(\G_t) = N \min_{\mathbf{1}^T\y =0, \y \neq 0} \frac{\sum_{m<n} w_{mn}(y_m-y_n)^2}{\sum_{m<n}(y_m-y_n)^2}. 
\end{align}
Here, the minimum is attained by the corresponding eigenvector $\breve{\y}$, that satisfies $\L_t\breve{\y} = a(\G_t)\breve{\y}$. Recall that $\E:=\{(m,n) \mid w_{mn} \in [\epsilon_w, 1]\}$, and observe that since $\G_t$ is connected, there exists a path $\mathcal{P}$ between any two nodes $m$ and $n$, such that
\begin{align}
(\breve{y}_m-\breve{y}_n)^2 &= [\sum_{(i,j) \in \mathcal{P}} \breve{y}_i - \breve{y}_j]^2 \! \leq (N-1)\!\!\sum_{(i,j) \in \mathcal{P}} \left(\breve{y}_i - \breve{y}_j\right)^2 \label{ineqmain} \\
& \leq (N-1)\sum_{(i,j) \in \E} \left(\breve{y}_i - \breve{y}_j\right)^2
\end{align} 
where, \eqref{ineqmain} holds since $\mathcal{P}$ may contain at most $N-1$ edges. Summing both sides over all edges in the graph, we have that
\begin{align}
\sum_{m <n} (\breve{y}_m - \breve{y}_n)^2 \leq    \frac{N(N-1)^2}{2}\sum_{(m,n) \in \E} \left(\breve{y}_m - \breve{y}_n\right)^2 \label{denombound}
\end{align}
Substituting \eqref{denombound} into \eqref{algdef} for $\y = \breve{\y}$, we have that
\begin{align}
a(\G) &= N\frac{\sum_{m<n} w_{mn}(\breve{y}_m-\breve{y}_n)^2}{\sum_{m<n}(\breve{y}_m-\breve{y}_n)^2} \\
&\hspace{-1cm}\geq \frac{2}{(N-1)^2} \frac{\sum_{(m,n)\in \E} w_{mn}(\breve{y}_m-\breve{y}_n)^2}{\sum_{(m,n)\in \E} (\breve{y}_m-\breve{y}_n)^2} \geq \frac{2\epsilon_w}{(N-1)^2} \label{algbound}
\end{align}
which is the required bound. If $\G_t$ is not connected, it holds for a component $\G_t^k$ with $p$ nodes that $a(\G_t^k) \geq 2\epsilon_w/(p-1)^2 \geq 2\epsilon_w/(N-1)^2$, so that we again have $a(\G_t) = \min a(\G_t^k) \geq 2\epsilon_w/(N-1)^2$, which is the desired result. 
\end{IEEEproof}

\vspace{-3mm}
\section{Proof of Lemma \ref{ktlem}}\label{lemproof}
Before proceeding with the proof, we state some basic results, and introduced necessary notation. In the subsequent analysis, we will repeatedly use the following inequalities \cite{golub2012matrix}
\begin{align} \label{frin}
\norm{\A\B} &\leq \norm{\A}_2\norm{\B} \leq \norm{\A}\norm{\B} 
\end{align}
where $\A$ and $\B$ matrices of compatible sizes. For notational brevity, $d_{mn} := \sqrt{\norm{\x_m-\x_n}^2+\epsilon_x}$ and $\check{d}_{mn} := \sqrt{\norm{\xc_m-\xc_n}^2+\epsilon_x}$, and note that $d_{mn}, \check{d}_{mn} \geq \sqrt{\epsilon}$. 

We begin by defining the \emph{total deviation} functions corresponding to $\K_t^1$ and $\K_t^2$ as
\begin{align}\label{sumdev}
\D^1_t(\X)&:= \sum_{\tau=1}^{t} \left(\L^{\dagger}_\tau\B^\epsilon_\tau(\X)\X - \Ex{\L^{\dagger}_\tau\B^\epsilon_\tau(\X)\X}\right) \\
\D^2_t(\X)&:= \sum_{\tau=1}^{t} \left(\L^{\dagger}_\tau\L_\tau-\Ex{\L^{\dagger}_\tau\L_\tau}\right)\X
\end{align}

The following lemma lists several preliminary results required in deriving the bounds in Lemma \ref{ktlem}. 

\begin{lem}\label{hovstatlem}
There exists $t_0 < \infty$, such that for all $t \geq t_0$, it holds that
\begin{subequations}
\begin{align}
\norm{\L_t^{\dagger}\B^\epsilon_t(\X)\X-\L_t^{\dagger}\B^\epsilon_t(\Xc)\Xc} &\leq C_3\norm{\X-\Xc}\label{bxlip} \\
\norm{\L_t^{\dagger}\B^\epsilon_t(\X)\X} &\leq C_4 \label{bxbound}\\
\norm{\J\Xh_t} &\leq C_5  \label{xbound}\\
\norm{\D^1_t(\X)} &\leq d^1_t \label{sumbound1}\\
\norm{\D^1_t(\X)-\D_t(\Xc)} &\leq \pi^1_t\norm{\X-\Xc}\label{sumlip1} \\
\norm{\D^2_t(\X)} &\leq d^2_t\label{sumbound2}\\
\norm{\D^2_t(\X)-\D^2_t(\Xc)} &\leq \pi^2_t\norm{\X-\Xc} \label{sumlip2}
\end{align}
\end{subequations}
where $\J = \I-\bm{11}^T/N$, $C_3$ and $C_4$ are constants, while the random variables $d^1_t$, $d^2_t$, $\pi^1_t$, and $\pi^2_t$ follow \eqref{dpit}. Results in \eqref{sumbound2} and \eqref{sumlip2} also require $\X$ to be such that $\norm{\J\X}\leq C_5$. 
\end{lem}
The proof organized into four steps, each considering one or more inequalities.

\begin{IEEEproof}[Proof of \eqref{bxlip} and \eqref{bxbound}]
Observe that the $m$-th row of $\B^\epsilon_t(\X)\X$ for each $t\geq 0$ can be written as
	\begin{align} \label{mrow}
	\left[\Bt(\X)\X\right]_{m,:} = \sum_{n\neq m} \frac{w_{mn}(t)\delta_{mn}(t)}{d_{mn}}\left(\x_m-\x_n\right) \nonumber
	\end{align}
	which implies that 
	\begin{align}
	\norm{\left[\Bt(\X)\X\right]_{m,:}} \leq \sum_{n\neq m} \abs{w_{mn}(t)\delta_{mn}(t)} .\leq NC_\delta.
	\end{align}
The bound in \eqref{bxbound} therefore follows from the use of \eqref{frin},
\begin{align}
\norm{\L_t^{\dagger} \Bt(\X)\X} &\leq \norm{\L_t^{\dagger}}_2\norm{\Bt(\X)\X}\leq \frac{N^2 C_{\delta}}{\epsilon_\L}. \label{ltbxbound}
	\end{align}
which yields $C_4 = N^2 C_\delta/{\epsilon_\L}$. Likewise, the $m$-th row of $\Bt(\X)\X-\Bt(\Xc)\Xc$ becomes
	\begin{align}
	&\left[\Bt(\X)\X-\Bt(\Xc)\Xc\right]_{m,:} \nonumber\\
	&= \sum_{n\neq m} w_{mn}(t)\delta_{mn}(t)\left(\frac{\x_m-\x_n}{d_{mn}}-\frac{\xc_m-\xc_n}{\check{d}_{mn}}\right).\label{bmrow}
	\end{align}
	Adding and subtracting the term $(\xc_m-\xc_n)/d_{mn}$ to each term within the summation in \eqref{bmrow}, it can be seen that 
	\begin{align}
	&\frac{\x_m-\x_n}{d_{mn}}-\frac{\xc_m-\xc_n}{\check{d}_{mn}} \nonumber\\
	&= \frac{\x_m-\xc_m}{d_{mn}} - \frac{\x_n-\xc_n}{d_{mn}} + (\xc_m-\xc_n)\left(\frac{1}{d_{mn}}-\frac{1}{\check{d}_{mn}}\right) \nonumber\\
	&= \frac{\x_m-\xc_m}{d_{mn}} - \frac{\x_n-\xc_n}{d_{mn}} + \frac{\xc_m-\xc_n}{\check{d}_{mn}}\frac{\check{d}^2_{mn}-d^2_{mn}}{d_{mn}(\check{d}_{mn}+d_{mn})}.
	\end{align}
Further, the term $\check{d}^2_{mn}-d^2_{mn}$ can be written compactly as
	\begin{align}
	\check{d}^2_{mn}-d^2_{mn} &= \xc_m^T\xc_m+\xc_n^T\xc_n - 2\xc_m^T\xc_n-\x_m^T\x_m-\x_n^T\x_n + 2\x_m^T\x_n \nonumber\\
	&\hspace{-1.5cm}= (\x_m -\x_n + \xc_m - \xc_n)^T(\x_m-\xc_m + \xc_n-\x_n) \label{xcmn}
	\end{align}
Consequently, it is possible to write \eqref{bmrow} as, 
\begin{align}
&\left[\Bt(\X)\X-\Bt(\Xc)\Xc\right]_{m,:} \nonumber\\
&=\sum_{n\neq m} w_{mn}(t)\delta_{mn}(t)\A_{mn} \left((\x_m-\xc_m) - (\x_n-\xc_n)\right) \nonumber
\end{align}
where the matrix $\A_{mn}$ is given by
\begin{align}\label{amn}
\A_{mn} &= \frac{1}{d_{mn}}\I + \frac{(\xc_m-\xc_n)(\x_m -\x_n + \xc_m - \xc_n)^T}{d_{mn}\check{d}_{mn}(\check{d}_{mn}+d_{mn})}. 
\end{align}
Thus, the full difference becomes 
\begin{align}
\Bt(\X)\X-\Bt(\Xc)\Xc = \A_t(\X,\Xc)\vect{\X-\Xc}
\end{align}
where the $(m,n)$-th $p \times p$ block of $\A_t(\X,\Xc)$ is given by
\begin{align}
\left[\A_t(\X,\Xc)\right] := \begin{cases} -\A_{mn}w_{mn}(t)\delta_{mn}(t) & m\neq n \\
\sum_{n\neq m} \A_{mn}w_{mn}(t)\delta_{mn}(t) & m = n
\end{cases}
\end{align}
Next, repeated use of the triangle inequality yields 
\begin{align}
\norm{\A_{mn}}^2 & \nonumber\\
&\hspace{-1cm}\leq \frac{2}{d^2_{mn}}\left(\norm{\I}^2 + \frac{\norm{\xc_m-\xc_n}^2}{\check{d}^2_{mn}}\frac{\norm{\xc_m-\xc_n + \x_m-\x_n}^2}{(\check{d}_{mn} + d_{mn})^2}\right) \nonumber
\end{align}
Here, it holds from the definition of $\check{d}_{mn}$ that $\norm{\xc_m-\xc_n}/\check{d}_{mn} \leq 1$. Similarly, it holds that 
\begin{align}
&\norm{\xc_m-\xc_n + \x_m-\x_n}^2 \\
&\leq \norm{\xc_m-\xc_n}^2 + \norm{\x_m-\x_n}^2 + 2\norm{\xc_m-\xc_n}\norm{\x_m-\x_n} \nonumber\\
&\leq \check{d}_{mn}^2 + d_{mn}^2 + 2\check{d}_{mn}d_{mn}=(\check{d}_{mn} + d_{mn})^2 
\end{align}
Therefore, the bound on $\norm{\A_{mn}}^2$ becomes 
\begin{align}
\norm{\A_{mn}}^2 &\leq \frac{2(N+1)}{\epsilon}
\end{align}
Similarly, it holds for $\norm{\A_t(\X,\Xc)}$ that
\begin{align}
\norm{\A_t(\X,\Xc)}^2 &\leq  C_{\delta}^2\sum_{m} \sum_{n\neq m}\norm{\A_{mn}}^2 + \left(\sum_{n\neq m}\norm{\A_{mn}}\right)^2 \nonumber\\
& \leq 3 C_\delta^2\sum_{m}\sum_{n\neq m} \norm{\A_{mn}}^2\\
& \leq 3 C_\delta^2\frac{N(N-1)(N+1)}{\epsilon} < C_\delta^2\frac{6N^3}{\epsilon_x}
\end{align}
which in turn, yields the bound 
\begin{align}
\norm{\A_t(\X,\Xc)}^2 \leq 6N^3\frac{C^2_\delta}{\epsilon_x}.
\end{align} 
The Lipschitz continuity of $\L_t^{\dagger}\B_t(\X)\X$ thus follows as 
\begin{align}
\norm{\L_t^{\dagger}\B_t^\epsilon(\X)\X - \L_t^{\dagger}\B_t^\epsilon(\Xc)\Xc} & \leq
\norm{\L_t^{\dagger}}_2\norm{\B_t^\epsilon(\X)\X - \B_t^\epsilon(\Xc)\Xc} \nonumber\\
&\hspace{-1cm}\leq \frac{NC_\delta}{\epsilon_\L} \sqrt{\frac{6N}{\epsilon_x}}\norm{\X-\Xc},
\end{align}
so that $C_3 = \frac{N C_\delta}{\epsilon_\L} \sqrt{\frac{6N}{\epsilon_x}}$.
\end{IEEEproof}

\begin{IEEEproof}[Proof of \eqref{xbound}]
Observe that $\L_t\J = \L_t$ and $\J\L_t^\dagger = \L_t^\dagger$. Right multiplying both sides of \eqref{smdsg} by $\J$, it follows that
\begin{align}
\J\Xh_{t+1} &= \J(\I - \mu \L_t^\dagger\L_t)\Xh_t + \mu \J\L_t^\dagger\B_t^\epsilon(\Xh_t)\Xh_t \\
&= (\J\J - \mu \J\L_t^\dagger\L_t\J)\Xh_t + \mu \L_t^\dagger\B_t^\epsilon(\Xh_t)\Xh_t \\
&= (\J - \mu \L_t^\dagger\L_t)\J\Xh_t + \mu \L_t^\dagger\B_t^\epsilon(\Xh_t)\Xh_t \label{vsmds} \\
&\hspace{-1cm}= (\J - \mu \L_t^\dagger\L_t)(\J-\mu \L_{t-1}^{\dagger}\L_{t-1})\J\Xh_{t-1}  + \mu \L_t^\dagger\B_t^\epsilon(\Xh_t)\Xh_t\nonumber\\
& + \mu(\J - \mu \L_t^\dagger\L_t)\L_{t-1}^\dagger\B_{t-1}(\Xh_{t-1})\Xh_{t-1} \label{vsmds2}
\end{align}
Continuing in a similar manner, taking norm on both sides of \eqref{vsmds2}, applying triangle inequality, and using \eqref{bxbound} yields
\begin{align}
&\norm{\J\Xh_{t+1}} \leq \norm{\Q_t^0}_2\norm{\J\Xh_0} + \mu(1+\sum_{\tau=1}^t\norm{\Q^\tau_t}_2)C_4 \label{vxineq1}
\end{align} 
where $\Q^\tau_t := \prod_{\kappa=\tau}^t(\J - \mu\L_t^\dagger\L_t)$. Next, from (\textbf{A4}), there exists some $t_0 < \infty$ and $\varrho < 1$ such that $\norm{\Q_{t}^\tau} \leq \varrho^{t-\tau+1}$ for all $t-\tau+1 \geq t_0$.  Since $\norm{\Q_t^\tau} \leq 1$ for all $t\geq \tau+1$, bound in \eqref{vxineq1} becomes
\begin{align}
\norm{\J\Xh_{t+1}} \leq C_x \varrho^{t} + \mu C_4 (1+t_0 + \frac{\varrho^t}{1-\varrho}) \nonumber\\
=C_x + \mu C_4(1+t_0 + \frac{1}{1-\varrho})=:C_5 
\end{align}
for all $t \geq t_0$. 
\end{IEEEproof}

\begin{IEEEproof}[Proof of \eqref{sumbound1} and \eqref{sumlip1}]
Observe that each term of $\D_t^1(\X)$ in \eqref{sumdev} is zero mean, and bounded as
\begin{align}
&\norm{(\L^{\dagger}_\tau\B^\epsilon_\tau(\X)\X - \Ex{\L^{\dagger}_\tau\B^\epsilon_\tau(\X)\X}} \\
&\leq \norm{(\L^{\dagger}_\tau\B^\epsilon_\tau(\X)\X} + \norm{\Ex{\L^{\dagger}_\tau\B^\epsilon_\tau(\X)\X}} \\
&\leq \norm{(\L^{\dagger}_\tau\B^\epsilon_\tau(\X)\X} + \Ex{\norm{\L^{\dagger}_\tau\B^\epsilon_\tau(\X)\X}} \leq 2C_4
\end{align}
The law of large numbers therefore implies that $\D^1_t(\X)/t \rightarrow 0$ almost surely as $t\rightarrow \infty$. This also implies that there exists $d^1_t$ such that $\norm{\D_t^1(\X)} \leq d_t^1$ and $d_t^1/t \rightarrow 0$ as $t \rightarrow \infty$.

The Lipschitz continuity of $\D_t^1(\X)$ can similarly be shown using \eqref{bxlip}. Towards this end, observe that
\begin{align}
\D^1_t(\X) - \D^1_t(\Xc) &=\sum_{\tau=1}^{t-1} \left(\L^{\dagger}_\tau\B^\epsilon_\tau(\X)\X - \Ex{\L^{\dagger}_\tau\B^\epsilon_\tau(\X)\X}\right) \nonumber\\
&\hspace{-0cm}- \sum_{\tau=1}^{t-1} \left(\L^{\dagger}_\tau\B^\epsilon_\tau(\Xc)\Xc - \Ex{\L^{\dagger}_\tau\B^\epsilon_\tau(\Xc)\Xc}\right) \nonumber\\
&\hspace{-1cm}= \sum_{\tau=0}^{t-1} \L^{\dagger}_\tau \left(\B^\epsilon_\tau(\X)\X - \B^\epsilon_\tau(\Xc)\Xc\right) \nonumber\\
&\hspace{-0cm} - \Ex{\L^{\dagger}_\tau \left(\B^\epsilon_\tau(\X)\X - \B^\epsilon_\tau(\Xc)\Xc\right)}\label{pibfdiff}
\end{align}
The vectorized version of the first term can be written as
\begin{align}\label{lbx}
&\vect{\L^{\dagger}_\tau\B^\epsilon_\tau(\X)\X - \L^{\dagger}_\tau\B^\epsilon_\tau(\Xc)\Xc} \nonumber\\
&\hspace{1cm}= \left(\I \otimes \L_\tau^{\dagger} \right)\vect{\B^\epsilon_\tau(\X)\X-\B^\epsilon_\tau(\Xc)\Xc} \\
&\hspace{1cm}= \left(\I \otimes \L_\tau^{\dagger} \right)\A_\tau(\X,\Xc)\vect{\X-\Xc}
\end{align}
Using a similar transformation on the second term of \eqref{pibfdiff}, the vectorized version of the right-hand side can be written as
\begin{align}
&\vect{\D_t^1(\X) - \D_t^1(\Xc)} \nonumber\\ &=\left(\sum_{\tau=0}^{t-1}\C_\tau(\X,\Xc)-\Ex{\C_\tau(\X,\Xc)}\right)\vect{\X-\Xc}
\end{align}
where $\C_\tau(\X,\Xc) = \left(\I \otimes \L_\tau^{\dagger} \right)\A_\tau(\X,\Xc)$ is bounded as $\norm{\C_\tau(\X,\Xc)} \leq \norm{\L_t^\dagger}_2\norm{\A_\tau(\X,\Xc)} \leq C_3$.
It is therefore possible to write
\begin{align}
\norm{\D^1_t(\X) - \D^1_t(\Xc)} \leq \pi_t \norm{\X-\Xc}
\end{align}
\begin{align}
\text{where,}\quad \pi_t  = \norm{\sum_{\tau=0}^{t-1} \C_\tau(\X,\Xc)-\Ex{\C_\tau(\X,\Xc)}}
\end{align}
Since the term within the norm is a bounded zero-mean random variable, it follows from law of large numbers that
\begin{align}
\frac{1}{t}\sum_{\tau=0}^{t-1} \C_\tau(\X,\Xc)-\Ex{\C_\tau(\X,\Xc)} \rightarrow \bs{0}
\end{align}
with probability 1 as $t\rightarrow \infty$. This also implies that $\pi_t/t \rightarrow 0$ almost surely as $t\rightarrow \infty$. 

\subsubsection{Proof of \eqref{sumbound2} and \eqref{sumlip2}} 
Observe that the zero mean random variable $\D_t^2(\X)$ can be written as 
\begin{align}
\D_t^2(\X) &= \sum_{\tau=0}^t \left(\L_\tau^\dagger\L_\tau - \Ex{\L_\tau^\dagger\L_\tau}\right)\J\X 
\end{align}
so that it follows form \eqref{xbound} that $\norm{\left(\L_\tau^\dagger\L_\tau - \Ex{\L_\tau^\dagger\L_\tau}\right)\J\X} \leq 2C_5$ for all $\X$ such that $\norm{\J\X}\leq C_5$. Invoking the law of large numbers as before, $\D_t^2(\X)/t \rightarrow 0$ almost surely as $t \rightarrow \infty$. Consequently, there exists $d_t^2$ such that $\norm{\D_t^2(\X)}\leq d_t^2$ and $d_t^2/t \rightarrow 0$ almost surely as $t \rightarrow \infty$.
 
In order to establish the Lipschitz continuity of $\D^2_t(\X)$, observe that $\D_t^2(\X) - \D_t^2(\Xc) = \C'_t(\X-\Xc)$, where
\begin{align}
\C'_t := \sum_{\tau=0}^t \L_\tau^\dagger\L_\tau - \Ex{\L_\tau^\dagger\L_\tau} \label{d2lipproof}
\end{align}
Since each summand in \eqref{d2lipproof} is zero mean and bounded, it holds from law of large numbers that  $\C'_t/t \rightarrow 0$ almost surely as $t \rightarrow \infty$. Consequently, there exists 
$\pi^2_t$ such that $\norm{\D_t^2(\X) - \D_t^2(\Xc)} \leq \pi^2_t\norm{\X-\Xc}$, and $\pi^2_t/t\rightarrow 0$ almost surely as $t \rightarrow \infty$. 
\end{IEEEproof}

\begin{IEEEproof}[Proof of Lemma \ref{ktlem}]
Bounds in \eqref{ktib} can be derived by observing that for $1\leq \tau\leq t$ and $\iota = 1, 2$, it holds that
\begin{align}\label{pibfd2}
\D^\iota_\tau(\Xh_{\tau}) - \D^\iota_{\tau-1}(\Xh_{\tau-1}) &= \K^\iota_\tau -\K^\iota_{\tau-1} + \nonumber\\
& \D^\iota_{\tau-1}(\Xh_{\tau}) - \D^\iota_{\tau-1}(\Xh_{\tau-1}).
\end{align}
Summing \eqref{pibfd2} over $\tau = 1, \ldots, t$, it follows that
\begin{align}
\D^\iota_t(\Xh_t) - \D^\iota_0(\Xh_0)  &= \K^\iota_t - \K^\iota_0 + \sum_{\tau=1}^{t} \left(\D^\iota_{\tau}(\Xh_{\tau+1}) - \D^\iota_{\tau}(\Xh_{\tau})\right) \nonumber 
\end{align}
Observing that $\K^\iota_0 = \D^\iota_0(\Xh_0)$, a bound on $\K^\iota_t$ can be derived by using \eqref{sumbound1} and \eqref{sumlip1} as follows:
\begin{align}
\norm{\K^\iota_t} &\leq \norm{\D^\iota_t(\Xh_t)} + \sum_{\tau=1}^{t} \norm{\D^\iota_{\tau}(\Xh_{\tau+1}) - \D^\iota_{\tau}(\Xh_{\tau})}\\
&\leq d^\iota_t + \sum_{\tau=1}^t \pi^\iota_\tau \norm{\Xh_{\tau+1}-\Xh_{\tau}} \\
&= d^\iota_t + \mu\sum_{\tau=1}^t \pi^\iota_\tau \norm{\L_{\tau}^{\dagger}\B^\epsilon_{\tau}(\Xh_{\tau})\Xh_{\tau} - \L_\tau^\dagger\L_\tau\Xh_{\tau}} \\
&\leq d^\iota_t + \mu\sum_{\tau=1}^t \pi^\iota_\tau \left(C_4 + \norm{\L_\tau^\dagger\L_{\tau}\J\Xh_{\tau}}\right) \\
&\leq d^\iota_t + \mu (C_4 + C_5) \sum_{\tau=1}^t \pi^\iota_\tau
\end{align}
so that $C_1 = C_2 = (C_4 + C_5)$ for $\iota = 1,2$. 

The bound on $\norm{\K^3_t}$ follows form applying triangle inequality on \eqref{kt3}, and using \eqref{bxlip} as follows:
\begin{align}
\norm{\K^3_t} &\leq \sum_{\tau=1}^{t-1}\norm{\Ex{\L_\tau^\dagger\B^\epsilon_\tau(\Xh_\tau)\Xh_\tau-\L_\tau^\dagger\B^\epsilon_\tau(\Xt_\tau)\Xt_\tau}} \\
	&\leq \sum_{\tau=1}^{t-1}{\Ex{\norm{\L_\tau^\dagger\B^\epsilon_\tau(\Xh_\tau)\Xh_\tau-\L_\tau^\dagger\B^\epsilon_\tau(\Xt_\tau)\Xt_\tau}}}\\
	&\leq \sum_{\tau=1}^{t-1}{C_3\norm{\Xh_\tau-\Xt_\tau}} =C_3\sum_{\tau=1}^{t-1}\norm{\Dl_\tau}
\end{align}

Finally, to show that $f_t(\mu) \leq f_T(\mu) \rightarrow 0$ for the interval $0 \leq t \leq T/\mu$, observe that for $\iota = 1, 2$, it holds that $\mu d^\iota_t \leq Td^\iota_t/t$. From \eqref{dpit}, it is known that given any $\varepsilon$, there exists $t_0(\varepsilon)$ and $C_d$ such that 
\begin{align}
\mathbb{P}\left[d^\iota_t/t \leq C_d\right]&=1 &\forall~t,\\
\text{and }~~~~\mathbb{P}\left[d^\iota_t/t \leq \varepsilon\right]&=1 &\forall~t>t_0(\varepsilon).
\end{align}
Such a $t_0(\varepsilon)$ exists within $[0,T/\mu]$ for all $\mu \leq T/t_0(\varepsilon)$. Therefore, given $\varepsilon$, if $t\leq t_0(\varepsilon)$, it holds that
\begin{align}
\mathbb{P}\left[\mu d^\iota_t \leq \varepsilon\right]=1 \label{pd}
\end{align}
 for all $\mu\leq \varepsilon / t_0(\varepsilon)C_d$. On the other hand, if $t>t_0(\varepsilon)$, \eqref{pd} holds for all $\mu\leq  T / t_0(\varepsilon/T)$. Combining the two cases, it holds that $ \max_{0\leq t\leq T/\mu} \mu d^\iota_t \rightarrow 0$, with probability one as $\mu \rightarrow 0$. 

For the other two terms, observe similarly that given $\varepsilon$, there exists $T_\varepsilon$ and $C_\pi$ such that 
\begin{align}
\mathbb{P}\left[\pi^\iota_t/t \leq C_\pi\right]&=1 &\forall~t,\\
\text{and }~~~~\mathbb{P}\left[\pi^\iota_t/t \leq \varepsilon\right]&=1 &\forall~t>T_\varepsilon.
\end{align}
Thus, given $\varepsilon$, if $t\leq T_\varepsilon$, it holds that
\begin{align}
\mathbb{P}\left[\mu^2\sum_{\tau=2}^t \pi^\iota_\tau \leq \varepsilon\right]&=1, \quad \forall  \mu, \quad \text{s.t}, \quad \mu \leq \frac{1}{T_\varepsilon}\sqrt{\frac{\varepsilon}{C_\pi}}. \label{pd2} 
\end{align}
Similarly, the result in \eqref{pd2} holds for $t>T_\varepsilon$ for all  $\mu \leq \frac{T}{T_{\varepsilon/T^2}}.$
\end{IEEEproof}

\section{Proof of Lemma \ref{avglem}}\label{avglemProof}
For notational convenience, let $\breve{\delta}_{mn}(t):=\frac{\delta_{mn}(t)}{\sqrt{\norm{\x_m-\x_n}^2+\epsilon_x}}$ and recall that $\bar{\delta}_{mn} = \Ex{\breve{\delta}_{mn}(t)}$. The proof is divided into two parts. In the first part, we consider the case when $\G_t$ is connected, so that $p = N$. In this case, the goal is to show that
\begin{align}
N\left[\Ex{\L_t^\dagger\B^\epsilon_t(\X)}\right]_{mn} & = \begin{cases} -\bar{\delta}_{mn} & m \neq n \\
\sum_{n \neq m} \bar{\delta}_{mn} & m = n.
\end{cases}
\end{align}
Since the graph is connected, it holds that $\L_t^\dagger = (\L_t + \b{11}^T/N)^{-1}-\b{11}^T/N$. Let $\psi_{mn}$ denote the $(m,n)$-th co-factor of $\L_t + \b{11}^T/N$ and $\Psi:=\text{det}(\L_t+\b{11}^T/N)$, so that $[\L_t^\dagger]_{mn}=\psi_{mn}/\Psi-1/N$. Since $\L_t^\dagger$ has zero row and column sums, we also have that $\sum_{n=1}^M\psi_{mn} = \Psi$. Therefore, expanding along the $m$-th row, the expression for $\Psi$ becomes
\begin{align}\label{Psi}
\Psi &= \sum_{n \neq m} w_{mn}(t)(\psi_{mm}-\psi_{mn}) + \frac{1}{N}\sum_{n=1}^M \psi_{mn} \\
&= \frac{N}{N-1}\sum_{n \neq m} w_{mn}(t)(\psi_{mm}-\psi_{mn})
\end{align} 
for each $1\leq m \leq N$. 
Straightforward manipulations allow us to conclude that
\begin{align}
\left[\L_t^\dagger\B^\epsilon_{t}(\X)\right]_{mn} \hspace{-0.3cm}= \frac{1}{\Psi}\begin{cases}  - \breve{\delta}_{mn}(t)w_{mn}(t)(\psi_{mm}-\psi_{mn}) & \hspace{-1.2cm} m\neq n \\
~~~- \sum_{k\neq m,n}w_{nk}(t)\breve{\delta}_{nk}(t)(\psi_{mn}-\psi_{mk}) \\
\sum\limits_{k\neq m}w_{mk}(t)\breve{\delta}_{mk}(t)(\psi_{mm}-\psi_{mk}) & \hspace{-1.2cm} m = n.
\end{cases}\nonumber
\end{align}

Next, we show that the random variables $\psi_{mn}$ and $\psi_{mk}$ are identically distributed for $n \neq k\neq m$. Without loss of generality, let $m = 1$. Also, let $\L_{i}^{nk}$ denote the $(N-2)\times(N-2)$ submatrix of $\L_t+\b{11}^T/N$ after the removal of rows $(1,i)$ and columns $(n,k)$. The Laplace expansion of $\psi_{1n}$ along the $k$-th column  yields 
\begin{align}
\psi_{1n} &= -\sum_{i \neq 1,n,k} (\frac{1}{N}-w_{ki}(t))(-1)^{n+i+k}\abs{\L_{i}^{nk}} \nonumber\\
&\hspace{-0.5cm}- (\frac{1}{N}-w_{kn}(t)) (-1)^{k}\abs{\L_{n}^{nk}} - (\frac{1}{N} + \sum_{i\neq k}w_{ki}(t))(-1)^{n}\abs{\L_{k}^{nk}} \nonumber\\
&= -\sum_{i \neq 1} (\frac{1}{N}-w_{ki}(t))(-1)^{n+i+k}\abs{\L_{i}^{nk}} \nonumber\\
&- (\sum_{i\neq k,n}w_{ki}(t)+2w_{kn}(t))(-1)^{n}\abs{\L_{k}^{nk}} \label{psin}
\end{align}
Likewise, the expansion of $\psi_{1k}$ along the $n$-th column  yields
\begin{align}
\psi_{1k} &= -\sum_{i \neq 1} (\frac{1}{N}-w_{ni}(t))(-1)^{n+i+k}\abs{\L_{i}^{nk}} \nonumber\\
&- (\sum_{i\neq k,n}w_{ki}(t)+2w_{kn}(t))(-1)^{k}\abs{\L_{n}^{nk}} \label{psik}
\end{align}
It can be seen that the first terms in \eqref{psin} and \eqref{psik} are identically distributed since $w_{ni}(t)$ and $w_{ki}(t)$ are identical (cf. (\textbf{A5})). Further, performing $n-k$ row exchanges on $\L_{n}^{nk}$, it is possible to obtain $\tilde{\L}_{n}^{nk}$ which only differs from $\L_{k}^{nk}$ in the $k$-th row. Indeed, the elements of the $k$-th row of $\tilde{\L}_{n}^{nk}$ are $\{(1/N-w_{ki}(t))\}_{i\neq k,n}$, while the elements of the $k$-th row of $\L_n^{nk}$ are $\{(1/N-w_{ni}(t))\}_{i\neq k,n}$. Since the determinant is linear in its rows, it follows that $\abs{\L_{n}^{nk}}$ and $\abs{\tilde{\L}_{n}^{nk}} = (-1)^{n+k}\L_{k}^{nk}$ are identically distributed. In summary, we have that the distributions of $\psi_{mn}$ and $\psi_{mk}$ are identical for all $k \neq n \neq m$. 

Next, define identical random variables $\chi_{mn} := w_{mn}(t)(\psi_{mm}-\psi_{mn})$ for each $n \neq m$, so that $\Psi = \frac{N}{N-1}\sum_{n\neq m} \chi_{mn}$. Since $\G_t$ is connected, it holds that $\Psi > 0$. Therefore from symmetry, we have that 
\begin{align} \label{chiexp}
\Ex{\frac{\chi_{mn}}{\Psi}} &= \frac{N-1}{N}\Ex{\frac{\chi_{mn}}{\sum_{n\neq m} \chi_{mn}}} = \frac{1}{N}
\end{align}
Further, using the fact that $\Ex{\chi_{mn}} = \Ex{\chi_{mk}}$ for each $k \neq n$, it can be seen that 
\begin{align}
\mathbb{E}\left[\L_t^\dagger\B^\epsilon_{t}(\X)\right]_{mn}= \frac{1}{N}\begin{cases}  - \bar{\delta}_{mn} & m\neq n \\
\sum\limits_{k\neq m} \bar{\delta}_{mk} &m = n.
\end{cases}\nonumber
\end{align}
which is the required result. 

Finally, if $\G_t$ consists of multiple connected components, the quantity $\L_t^\dagger\B^\epsilon_{t}(\X)$ is a permuted version of the block-diagonal matrix with $N/p$ block matrices of size $p \times p$ each. Let $\Psi^j$ denote the determinant of $j$-th block, and the random variables $\chi^j_{mn}$ be similarly defined block-wise. Proceeding along similar lines, it can be seen that 
\begin{align} \label{chiexp2}
\Ex{\frac{\chi^j_{mn}}{\Psi^j}} &= \frac{p-1}{p}\Ex{\frac{\chi^j_{mn}}{\sum_{n\neq m} \chi^j_{mn}}} = \frac{1}{p}.
\end{align}
Consequently, $[\L_t^\dagger\B^\epsilon_{t}(\X)]_{mn}$ is non-zero if and only if the node pair $(m,n)$ belong to the same component, and is zero otherwise. From $(\textbf{A5})$, we have that the probability that a given pair of nodes $(m,n)$ belongs to the same connected component is given by $(p-1)/(N-1)$, yielding the required expression
\begin{align}
\mathbb{E}\left[\L_t^\dagger\B^\epsilon_{t}(\X)\right]_{mn}= \frac{p-1}{p(N-1)}\begin{cases}  - \bar{\delta}_{mn} & m\neq n \\
\sum\limits_{k\neq m} \bar{\delta}_{mk} &m = n.
\end{cases}\nonumber
\end{align}


\bibliographystyle{IEEEtran}
\bibliography{IEEEabrv,references}

\end{document}